\documentclass[12pt,oneside,a4paper,leqno]{article}
\usepackage{psfrag}
\usepackage{subfigure}

\usepackage[cmtip,arrow,line,all]{xy}
\CompileMatrices
\usepackage{geometry}
\usepackage{graphicx}
\usepackage{amsfonts,amssymb}
\usepackage[centertags]{amsmath}
\usepackage{amsthm}

\newcounter{theorems}

\theoremstyle{plain}
\newtheorem{theoAA}[theorems]{Theorem}

\swapnumbers
\newcounter{lemma}
\numberwithin{equation}{section}

\newtheoremstyle{par}%
     {\topsep}%
     {\topsep}%
     {\itshape}%
     {}%
     {\bfseries}%
     {}%
     {.5em}%
     {}%
\theoremstyle{plain}
\numberwithin{equation}{section}

\newtheorem{propo}[equation]{Proposition}

\theoremstyle{definition}
\newtheorem{defi}[equation]{Definition}
\newtheorem{example}[equation]{Example}

\theoremstyle{remark}
\newtheorem{remark}[equation]{Remark}

\theoremstyle{par}
\newtheorem{nr}[equation]{}
\newtheorem{lemma}[equation]{}
\makeatletter
\def\tagform@#1{\maketag@@@{\ignorespaces#1\unskip\@@italiccorr}}

\makeatother
\newcommand{\minus}{\smallsetminus}
\newcommand{\CC}{\mathbb{C}}
\newcommand{\R}{\mathbb{R}}
\newcommand{\ze}{\mathbb{Z}}
\newcommand{\from}{\colon}
\newcommand{\X}{\mathcal{X}}
\newcommand{\T}{\mathbb{T}}
\newcommand{\I}{\mathbb{I}}
\newcommand{\action}{\mathcal{A}}

\newcommand{\n}{\mathbf{n}}
\newcommand{\nn}{\mathbf{n}}
\newcommand{\kk}{\mathbf{k}}
\theoremstyle{par}

\newtheorem{prop}[equation]{}

\theoremstyle{plain}

\newtheorem{theorem}[equation]{Theorem}
\newtheorem{corollary}[equation]{Corollary}

\theoremstyle{definition}
\newtheorem{definition}[equation]{Definition}

\newcommand{\N}{\mathbb{N}}
\begin{document}
\pagenumbering{arabic}
\title{%
Symmetry groups of the planar $3$-body problem and 
action--minimizing trajectories}

\author{Vivina Barutello%
\footnotemark[1]\footnote{%
Dipartimento di Matematica e Applicazioni, 
Universit\`a di Milano Bicocca, 
via R. Cozzi, 53;
20125 Milano, Italy. \
email: \emph{vivina.barutello@unimib.it}, 
\emph{davide.ferrario@unimib.it},
\emph{susanna.terracini@unimib.it}},
\  Davide L.~Ferrario%
\footnotemark[1]%
\ \ and Susanna Terracini\footnotemark[1]%
}
\date{\today}
\maketitle

\begin{abstract}
We consider periodic and quasi-periodic solutions
of the three-body problem with homogeneous potential
from the point of view of equivariant calculus of variations.
First, we show that symmetry groups of the Lagrangian
action functional can be reduced to groups in a finite explicitly given
list, after a suitable change of coordinates. Then,
we show that local symmetric minimizers  are 
always collisionless, without any assumption on the group 
other than the fact that collisions
are not forced by the group itself. 
Moreover, we describe some properties of the resulting 
symmetric collisionless minimizers 
(Lagrange, Euler, Hill-type orbits and Chenciner--Montgomery 
figure-eight).
 
\vspace{12pt}
\noindent
\emph{MSC Subj. Class}: 
Primary 
70F10 (Mechanics of particles and systems: $n$-body problems);
Secondary 
70G75  (Mechanics of particles and systems: Variational methods).
37C80  (Dynamical systems and ergodic theory: Symmetries, equivariant dynamical systems),
70F16 (Mechanics of particles and systems: 
Collisions in celestial mechanics, regularization),

\vspace{12pt}
\end{abstract}
\section{Introduction and main results}
\label{sec:intro}

Among all periodic solutions of the planar $3$-body problem,
relative equilibrium motions -- the equilateral Lagrange
and the collinear Euler-Moulton solutions --
are definitely the simplest and most 
known. They are both endowed with an evident 
symmetry ($SO(2)$ and $O(2)$ respectively), 
that is, they are equivariant with respect to the symmetry group 
of dimension $1$ acting as $SO(2)$ (resp.\ $O(2)$)
on the time circle and on the 
plane, and trivially on the set of indexes $\{1,2,3\}$. 
In fact, they are minimizers of the Lagrangian
action functional in the space of all loops having 
their same symmetry group.
Hence,  $G$-equivariant minimizers
for the action functional (given a symmetry group $G$)
can be thought as the natural generalization of relative
equilibrium motions. 
The purpose of the paper  is to apply systematically 
the $G$-equivariant calculus of variations 
to the planar $3$-body problem -- 
with homogeneous interaction potential.
We will develop a theory and 
describe all the main aspects involved in
the $G$-equivariant minimization:
in particular,
the 
classification of  all possible symmetry
groups,  the 
description of the main properties of their minimizers 
and the proof of the fact that $3$-body minimizers are always 
collisionless.

In the past few years  some variational methods have been exploited 
by several authors in the search
of new periodic solutions as symmetric minimizers  for the 
$n$-body problem
(see e.g. 
\cite{%
monchen,
FT2003,
marchal,
montgomery,
mont_prepr}
for some recent results and references,
\cite{%
amco,%
bruno,%
chenICM,%
diacu,%
hsiang2004,%
marchal_book,%
meyer,
montgomery_contmath,%
sbano%
}
for more details on the three-body problem or  variational methods,
and 
\cite{%
amco94,
bara91,
bessi,%
cotizelati,
MR95k:58135,
hsiang1994,
serter1,
serter2}
for a similar approach dating back to 
the early nineties).
A major problem in the equivariant minimization  is 
that minimizers might be  \emph{a priori} colliding trajectories.
To exclude this possibility several arguments have been introduced
and used in recent literature. 
In \cite{FT2003} the authors proposed a class of symmetry groups with the
property that all local minimizers are collisionless (groups 
with the \emph{rotating circle property}), for the general
$\alpha$-homogeneous $n$-body problem in dimension $d\geq 2$.
The key-step is a generalization of the averaging technique
introduced by Marchal and exposed in 
\cite{chenICM}.
In spite of its generality, the main theorem of \cite{FT2003}
cannot be applied to some relevant symmetry groups, 
such as the symmetry group of the Chenciner--Montgomery
eight-shape orbit \cite{monchen}.
 
Besides that of collisions, the major problems 
occurring in the variational approach are the following:
first, the minimum has to be achieved (this requires the condition
of \emph{coercivity} of the action functional on the space
of symmetric loops; see~\ref{existence} below). 
This depends on the group $G$ 
and possibly on the angular velocity
$\omega$ of the rotating frame (when we consider 
these coordinates). Roughly speaking, 
a group for which the action functional
is not coercive for all angular velocities 
will be termed \emph{fully uncoercive}
(see~\ref{def:fullyunc} below for the precise definition). 

As we already pointed out, one has 
to prove that the minimizer is
collisionless; in particular, this obviously excludes 
the case when collisions
are forced by symmetries.  
A group will be termed \emph{bound to collisions}
if  every equivariant loop has a collision (see~\ref{BTC} below).  Finally, it
makes sense to investigate whether the minimizer is necessarily a homographic
solution, and therefore we would like to exclude those groups for which all
equivariant loops are rotating configurations, as $SO(2)$ and $O(2)$ above.
Such groups will be termed \emph{homographic} (see~\ref{homG} below).

Since fully uncoercive, bound to collisions or homographic groups
are unsuitable to this minimization approach, we are interested
in the classification of all the other groups.
Our main results are the following theorems:
\begin{theoAA}
\label{MT1}
Let $G$ be a symmetry group of the Lagrangian action functional in the
$3$-body problem. Then
$G$ is either bound
to collisions, fully uncoercive, homographic, or,
up to a change of rotating frame, 
it is conjugated to one of the
symmetry groups listed in table~\ref{table:class} (RCP stands for Rotating
Circle Property and HGM for Homographic Global Minimizer).
\end{theoAA}

\begin{table}[ht]
\caption{Planar symmetry groups with trivial core}
\label{table:class}
\begin{center}
\begin{tabular}{l|c|c|c|c|c|c}
\hline
\emph{Name} & $|G|$ & \emph{type R} & \emph{action type} &  \emph{trans.\ dec.}  & RCP & HGM \\
\hline
Trivial              & $1$ & yes &  & $1+1+1$ & yes & yes\\
Line                 & $2$ & yes & brake  & $1+1+1$ & (no) & no \\
$2$-$1$-choreography & $2$ & yes & cyclic & $2+1$  & yes & no\\
Isosceles            & $2$ & yes & brake &  $2+1$ & no & yes\\
Hill                 & $4$ & yes & dihedral &  $2+1$ & no & no\\
$3$-choreography     & $3$ & yes & cyclic & $3$ & yes & yes\\
Lagrange             & $6$ & yes & dihedral & $3$ & no & yes\\
$C_6$                & $6$ & no & cyclic & $3$ & yes & no\\
$D_6$                & $6$ & no & dihedral & $3$ & yes & no\\
$D_{12}$             & $12$ & no & dihedral & $3$ & no & no\\
\end{tabular}
\end{center}
\end{table}

\begin{figure}[ht]
\caption{The poset of symmetry groups for the planar $3$-body problem}
\label{fig:poset}
\xymatrix{%
& & &  {D_{12}}   &  & {12} \\
& & {\mbox{Lagrange}} 
\ar@/^/[ur]  
&  {C_6} 
\ar@{->}[u]  
 &   {D_6} 
\ar@/_/[ul]  
& {6} \\
& {\mbox{Euler--Hill}} & & & & {4}\\
&  & & {\mbox{Choreography}} 
\ar@/^/[uul]  
\ar@{->}[uu] 
\ar@/_/[uur]
& & {3}\\
{\mbox{Line}}
\ar@/^/[uur]  
 & {\mbox{$2$-$1$-choreography}} 
\ar@{->}[uu]  
& {\mbox{Isosceles}}
\ar@/_/[uul]  
\ar@{->}[uuu]
& & & {2} \\
& &  
{\mbox{Trivial}}
\ar@/^/[ul]
\ar@/^/[ull] 
\ar@{->}[u]
\ar@/_/[uur]
&  &  & {1} \\
}
\end{figure}
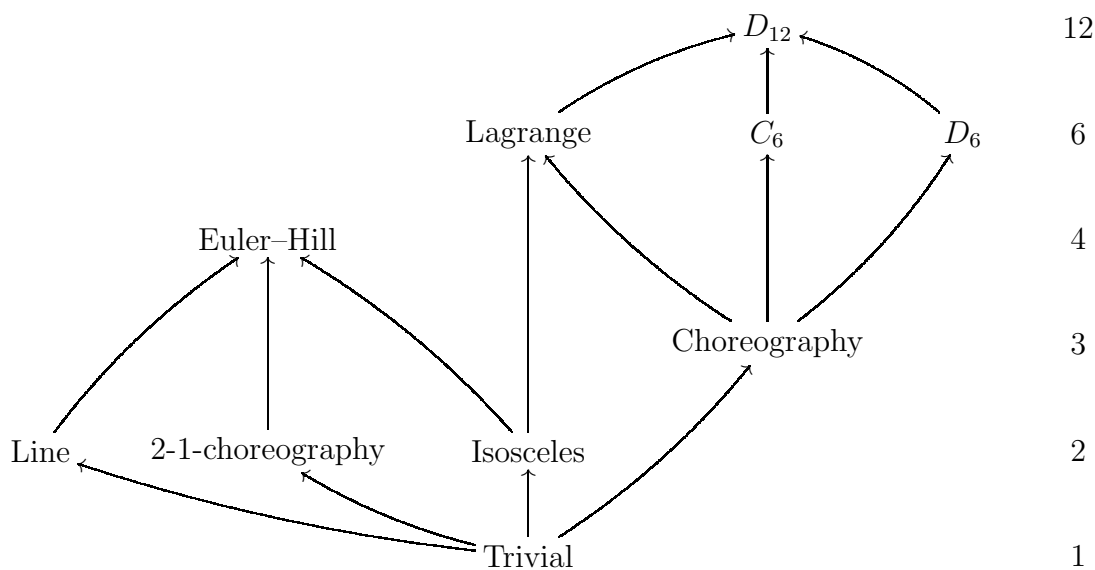

\begin{theoAA}
\label{MT2}
Let $G$ a symmetry group of the Lagrangian in the $3$-body
problem (in a rotating frame or not). If $G$ is not bound to collisions,
then any local minimizer is collisionless.
\end{theoAA}

The sentence ``up to  a change of rotating
frame'' in Theorem~\ref{MT1} requires a few words of explanation.
Consider a symmetry group $G$ for the Lagrangian functional:
if, for every angular velocity, 
$G$ is a symmetry group for the Lagrangian functional in the rotating frame,
then we will say that $G$ is \emph{of type R} (see the equivalent
definition~\ref{def:typeR} below).
This is a fundamental property for symmetry groups. In fact,
if $G$ is \emph{not} of type R, it turns out (see~\ref{lemma:nottypeR}) that 
the angular momentum of all $G$-equivariant trajectories vanishes (and thus,
that no homographic solution can occur).
For example, the dihedral symmetry group used by Chenciner and Montgomery
in \cite{monchen} is not of type R and thus, as already remarked
in their paper, the figure-eight  solution has zero angular momentum.
Only three groups are not of type R, for the planar $3$-body problem.
We will give a detailed description of each of them in section~\ref{sym:eights}, together 
with their minimizers; in this way we will 
partially answer to the open question (posed by Chenciner)
whether their minimizers coincide or not: for two of them ($D_6$ and $D_{12}$)
the answer is yes (under certain conditions, see~\ref{prop:D6D12} below).

On the other hand, if $G$ is of type R, then 
a suitable choice of the angular velocity of the rotating frame allows
us to consider only the minimal groups, listed in table~\ref{table:class}. 
Furthermore, given such a minimal group $G$ and a rational $\omega$,
$G$-symmetric trajectories (minimizers in the $\omega$-rotating  frame) yield
$G'$-symmetric trajectories for a different (generally bigger) group $G'$, which are  minimizers 
in the inertial frame (see remark~\ref{rem:typeR}).
Of course, when the group is of type R a key point 
is to understand whether the minimizer
is homographic or not. This indeed happens  every time the rotating Lagrange
configuration is $G$-equivariant. In the 
opposite cases,
whenever the Lagrange motion is \emph{not} $G$-equivariant, it is always
possible to find suitable choices of masses $m_i$ and angular velocity
$\omega$ such that the minimizer is not homographic; it turns out that these
are Hill-like solutions.
We emphasize that  our approach to exclude collisions is purely local,
as in \cite{serter1,FT2003}.
We do not 
need action estimates on colliding trajectories: in contrast we shall exhibit
local variations around parabolic ejection-collision solutions (which, as
shown in \cite{FT2003}, are the blow-ups of possible colliding minimizers).
This approach brings two advantages: first, it allows the extension of the
result in \cite{FT2003}
to a unified treatment of all symmetry groups, and also it works with all
$\alpha$-homogeneous potentials -- actually, this approach can
be extended to a larger class of potentials, under only
reasonable local assumptions at singularities.  Moreover, it 
allows one to 
prove existence of collisionless minimizers without any 
level estimate on the minimal actions of colliding trajectories,
that in the literature have been often obtained numerically (see 
for example 
\cite{montgomery_contmath,chenICM,monchen,chen2003_I,chen2003_II}).
On the other hand, our local variations are not suitable to be 
used in a context where constraints are of a different nature,
such as homotopy/homology/topology constraints \cite{montgomery}.

\section{Preliminaries}
\label{BFT:s2}

\subsection{Settings and notation} 
Let $E=\R^2\cong {\mathbb C}$ the
2-dimensional Euclidean space and $0 \in \R^2$ its origin. Let $m_1,m_2,m_3$
be positive real numbers and $\X$ the configuration space of three point
particles with masses $m_i$ respectively with center of mass in 0, then
\begin{equation} \label{configurationspace} 
\X =  \left\{
x=(x_1,x_2,x_3)\in E^3 : \sum_{i=1}^3 m_i x_i =0 \right\}.  
\end{equation} 
We will denote by $\Delta_{i,j}=  \left\{ x \in \X : x_i = x_j \right\}$
the collision set of the $i$-th and the $j$-th particle and with $\Delta =
\bigcup_{i,j}\Delta_{i,j}$, the \emph{collision set} in $\X$, i.e.
\begin{equation}
\label{collision:set}
\Delta =  \left\{ x \in \X : 
\exists i,j  \mbox{ \ such that \  }  x_i = x_j \right\}.
\end{equation}

Let ${\bf k}$ be a subset of the index set ${\bf n}=\{1,2,3\}$ and 
let ${\bf k'}$
denote its complement in  ${\bf n}$, then $\Delta_{{\bf k},{\bf
k'}}=\bigcup_{i\in{\bf k},j\in{\bf k'}}\Delta_{i,j}$ and $\Delta_{{\bf k}} =
\bigcap_{i,j \in{\bf k}}\Delta_{i,j}$.

Let $\alpha$ be a given positive real number. We consider the 
\emph{potential function} 
(opposite to the potential energy) and the \emph{kinetic energy}
defined respectively on $\X$ and on the tangent bundle of $\X$ as
\begin{eqnarray}
\label{pot:function}
U(x)&=& \frac{m_1 m_2}{|x_1-x_2|^\alpha}+
\frac{m_1 m_3}{|x_1-x_3|^\alpha}+\frac{m_2 m_3}{|x_2-x_3|^\alpha},\\
\label{kin:energy}
K(x,\dot{x})&=&
\frac{1}{2}\left( m_1|\dot{x}_1|^2 + 
m_2|\dot{x}_2|^2 + m_3|\dot{x}_3|^2 \right).
\end{eqnarray}

Associated to $U$ and $K$ there are the \emph{Newton equations}
\[
\label{eq:neq}
m_i \ddot{x}_i
=\frac{\partial U}{\partial x_i}. 
\]
When we work in a uniform rotating frame the kinetic energy needs to be
changed into the corresponding form
\begin{equation}
\label{kin:energy:rot}
K_\omega(x,\dot{x})= \sum_{i=1}^{3}\frac{1}{2} m_i|\dot{x}_i -J\omega x_i|^2.
\end{equation}

where $J$ is the complex unit and $\omega$ is the angular velocity. We suppose
that the origin (which coincide with the center of mass) is the point of the
plane fixed by the rotation. The Lagrangian is 
\begin{equation}
\label{lagrangian}
L_\omega(x,\dot{x})=L_\omega=K_\omega+U.
\end{equation}

Let $\T \subset \R^2$ denote a circle in $\R^2$ of length
$T=|\T|$. Moreover, let $\Lambda = H^1(\T,\X)$ be the
Sobolev space of $L^2$ loops $\T\rightarrow\X$ with $L^2$
derivative. It is an Hilbert space with the scalar product 
\[
x \cdot y = \int_{\T}(x(t)y(t)+\dot{x}(t)\dot{y}(t))dt
\]
and the corresponding norm denoted by $\|\cdot\|$. An equivalent norm is given
by 
\[
\|x\|' = \left( \sum_{i=1}^3 
\left(\int_{\T}\dot{x}^2_i(t)dt + [x_i]^2 \right)\right)^{\frac{1}{2}},
\]
where $[x]$ is the average $[x]=\int_{\T}x(t)dt$. 
For every $x=x(t) \in \Lambda$ we denote
$x^{-1}\Delta \subset \T$ the set of \emph{collision times}. The point
$t_0 \in \T$ is an 
\emph{isolated collision} for $x \in \Lambda$ when
$t_0 \in x^{-1}\Delta$ and there exists a neighborhood of $t_0$,
$(t_0-\epsilon, t_0+\epsilon)$, such that $(t_0-\epsilon, t_0+\epsilon) \cap
x^{-1}\Delta = \{t_0\}$. 

A \emph{colliding cluster} of particles at $t=t_0$ is a subset ${\bf k}
\subseteq \{1,2,3\}$ such that $x(t_0) \in \Delta_{{\bf k}}$ and $x(t_0)
\notin \Delta_{{\bf k},{\bf k'}}$.

The \emph{action functional} is the positive-definite function 
$\action_\omega: \Lambda \rightarrow \R \cup {\infty}$ defined by
\begin{equation}
\label{action:functional}
\action_\omega(x)= 
\int_{\T}L_\omega(x(t),\dot{x}(t))dt, \quad \omega \in \R. 
\end{equation}
The action functional $\action_\omega$
is of class $C^1$ on the subspace of collision-free
loops of $\Lambda$; hence collisionless critical points of 
$\action_\omega$ in
$\Lambda$ are $T$-periodic $C^2$-solutions of the 
associated Euler--Lagrange equations
\begin{equation}
\label{newton:eq}
m_i \ddot{x}_i
= 2J\omega \dot x_i+\omega^2 x_i +
\frac{\partial U}{\partial x_i};
\end{equation}
these solutions will be called \emph{classical solutions}, 
and, in the non-rotating frame, they satisfy Newton equations~\ref{eq:neq}.

\begin{remark}
In~\ref{configurationspace} we assumed that configurations
have center of mass in $0$. This constraint requires 
a few words of comment. Indeed, while 
the action functional $\action_0$  is invariant under
translations, the same is not true for $\action_\omega$
with $\omega\neq 0$. This means that the center of mass constraint
is in general not a natural  reduction. However,
the convexity property of the kinetic quadratic form
with respect to the center of mass
implies that all local minima of $\action_\omega$ 
have center of mass in $0$.
\end{remark}

\subsection{Symmetry constraints}
\label{sec:symmcon}
Let $G$ be a finite group, acting on a space $X$, i.e.\ there 
exists a map $G
\times X \rightarrow X$, $(g,x)\mapsto gx$, such that $(g_1,g_2 x) \mapsto
(g_2 \cdot g_1)x$, where $\cdot$ is the internal operations of $G$. The space
$X$ is then called \emph{$G$-equivariant}. Let $x \in X$, the \emph{fixer}
or
the \emph{isotropy group} of $x$ in $G$ is the set $G_x = \{ g \in G : gx =
x\}$. When $H$ is a subgroup of $G$, then the space $X^H \subset X$ consists
of all points $x \in X$ that are fixed by $H$, 
$X^H = \{ x \in X : G_x \supset
H\}$.

Consider a finite group $G$ and the 2-dimensional orthogonal representations
$\tau, \rho: G \rightarrow O(2)$ of $G$. By $\tau$ and $\rho$, $G$ acts on the
time circle $\T\subset {\mathbb R}^2$ and on the Euclidean space $E$,
respectively. Moreover, we can consider the group homomorphism $\sigma : G
\rightarrow \Sigma_3$, where $\Sigma_3$ is the group of all permutations of
$3$ elements, and we endow the set of indexes ${\bf n}= \{ 1,2,3\}$ with a $G$
action. The homomorphism $\sigma$ satisfies the property
\begin{equation}
\label{equalm}
\forall g \in G : \left( \sigma(g)(i)=j \Rightarrow m_i=m_j \right).
\end{equation}
Given $\rho$ and $\sigma$ with property~\ref{equalm}, $G$ acts orthogonally on the
configuration space $\X$ by
\begin{equation}
\label{G-action}
\forall g \in G, \quad g \cdot (x_1,x_2,x_3) = 
\left( \rho(g)x_{\sigma(g^{-1})(1)},
\rho(g)x_{\sigma(g^{-1})(2)},
\rho(g)x_{\sigma(g^{-1})(3)}   
 \right).
\end{equation}
for every $g \in G$. Understanding the action of $G$ as a space linear
transformation or a permutation on the index set, we can write
\begin{equation}
\label{G-action_short}
\forall i \in {\bf n}, \quad (gx)_i = gx_{g^{-1}i}.
\end{equation}
Furthermore, we can use the representation $\tau$ and the action of G on
$\X$ given by~\ref{G-action}, to define the action of $G$ on the set of
loops $\Lambda$ as follows
\begin{equation}
\label{G-actionL}
\forall g \in G, \,\, \forall t \in \T, \,\, 
\forall x \in \X, \quad (g \cdot x)(t) = gx(g^{-1}t).
\end{equation}
The loops in $\Lambda^G$ are called \emph{equivariant loops} and are the loops
$\T \rightarrow \X$ fixed by $G$; $\Lambda^G$ is a closed
linear subspace of $\Lambda$. 

\begin{defi}
Consider the action of $G$ on the index set $\n=\{1,2,3\}$.
The \emph{transitive decomposition} of $\n$ 
is the decomposition of $\n$ into its $G$-orbits
or, with an abuse of notation, their lengths. For example,
if $G$ acts transitively on $\n$ then $\n$ is decomposed into
a unique orbit of length $3$. If $G$ acts trivially, then the
decomposition is of type $1+1+1$, since $\n$ is decomposed
into its subsets of length $1$.
\end{defi}

In the sequel, if not explicitely stated,
we will assume the following additional hypotheses 
on the group action:
\begin{enumerate}
\item 
\label{hh1}
$\ker \tau \cap \ker \rho \cap \ker \sigma = 1$;
\item There exists no proper linear subspace $E' \subsetneq E$ such that 
\label{hh2}
\[
\forall i \in {\bf n}, \forall x \in \Lambda^G, \forall t \in \T x_i(t) \in E' \subsetneq E;
\]
\item 
\label{hh3}
There is not an integer $k \neq \pm 1$ such that
\[
\forall x \in \Lambda^G, \exists y \in \Lambda : 
\forall t \in \T,  x(t)=y(kt).
\]
\end{enumerate}
In fact, it is easy to show that there is a canonical
isomorphism $\Lambda^G \cong \Lambda {\bar G}$
where $\bar G = \frac{G}{\ker \tau \cap \ker \rho \cap \ker \sigma}$,
and so without loss of generality~\ref{hh1} can be assumed to be true.
Furthermore, if~\ref{hh2} holds, then it is possible 
to consider the $n$-body problem on the subspace
$E'\subset E$. Finally, condition~\ref{hh3} simply means 
that all loops are periodic of a period that is sub-multiple of $2\pi$,
and therefore up to rescaling time it is possible
to define another (smaller) group $G$. We will see later 
in section~\ref{sec:typeR} how to use this idea
to reduce the number of groups in the classification 
(condition~\ref{hh3} is actually equivalent to 
definition~\ref{defi:redundant} below). 
\begin{definition}
\label{homG}
A symmetry group is termed \emph{homographic} if all equivariant loops 
are constant up to similarity.
\end{definition}

\begin{example}
Consider a homographic group $G$ and a loop $x$ similar to  a  certain
configuration $\xi$. Without loss of generality we can assume that 
$\sum_i m_i \xi_i^2= 1$. Thus, elements of $\Lambda$ have the form
$x(t) = z(t) \xi$, 
for a complex valued function  $z(t)$ and a configuration $\xi$;
the 
action functional $\action$ is actually a function only of $z$ and $\xi$.
It is easy to see that if $x(t) = z(t) \xi$  is a critical point of 
$\action$ then $z(t)$ solves the Kepler problem and $\xi$ is a central
configuration.
We recall here the definition of central configuration (further details
can be found in \cite{Mo90}).

\begin{defi}\label{defi:cc}
A configuration $\xi$ is said a \emph{central configuration}
if $\xi$ is a critical point for the potential $U$ on the 
ellipsoid $\sum_i m_i \xi_i^2= 1$. 
\end{defi}

Central configurations for the $3$-body problem are well-known.  Up to
similarity there is only one non-collinear central configuration, termed the
\emph{Lagrange} configuration (it is an equilateral triangle), which minimizes
$U(x)$ in the ellipsoid $\sum_i m_i \xi_i^2 = 1$.  On the other hand, there
are $3$ collinear central configurations, termed \emph{Euler-Moulton}
configurations (or simply Euler configurations), which are saddle points for
$U$.

Associated with a planar 
central configuration there is a relative equilibrium
motion, as we already pointed out in the introduction.  Now we prove that they
are equivariant critical points for the action.  Consider first the group
$G=SO(2)$ acting canonically on the time circle $\T$ and the Euclidean plane
$\R^2\cong \mathbb{C}$, and trivially on the index set $\{1,2,3\}$.  The
equivariant constraints are the following: for every $\theta \in \R$ and every
$t\in \R$ one has $x(t+\theta) = e^{J\theta} x(t)$.  Thus the group is
homographic and one can easily verify that the action is an affine function of
$U(\xi)$; therefore critical points of $\action^G$ are the rotating central
configurations above (for every choice of masses). The minimal action is
attained on the Lagrange motion (since the Lagrange configuration minimizes
the potential $U$).

Furthermore, consider the group $G=O(2)$ acting again canonically.  In this
case the equivariant constraints are as follows: for every $\theta \in \R$ and
every $t\in \R$ one has $x(t+\theta) = e^{J\theta} x(t)$, and, additionally,
it must be that $x(-t) = \overline x(t)$ for every $t\in \R$ (where $\overline
x$ is the complex conjugate of $x \in \mathbb{C}$).  The latter constraint
implies that  $x(0)$ belongs to the real line, that is, the rotating
configuration must be collinear. Hence, the only critical points are rotating
Euler-Moulton configurations.  In the paper we will deal with finite groups of
symmetries since, as shown above, infinite symmetry groups are always
homographic.  
\end{example}

From our point of view homographic groups are not interesting, since 
they yield critical points which are already well-known. Another class
of groups that are unsuitable to our approach are groups
for which symmetry constraints force the occurrence of collisions.

\begin{definition} \label{BTC} 
Let $G$ be a finite group acting on the set $\Lambda$,
we say that $G$ is \emph{bound to collisions}
if every $G$-equivariant loop has at least a collision time
\[
\forall x \in \Lambda^G, 
x^{-1}\Delta \neq \emptyset
\]
\end{definition}

\begin{example}
Consider the group $G$ of order $2$ acting on $\T$ 
as a reflection along a line, 
trivially on $\R^2$ and by permuting $1$ and $2$ in 
the index set $\{1,2,3\}$. Thus the equivariant constraints are:
\[
x_1(-t) = x_2(t), \ \ 
x_2(-t) = x_1(t), \ \
x_3(-t) = x_3(t).
\]
This implies that, for $t=0$ or $t=\pi$ one has
$x_1(t)=x_2(t)$, and hence all equivariant loops have collisions at times $0$
and $\pi$.
\end{example}

\begin{example}
Consider the dihedral group $K=D_6$ of order $6$ 
with generators $g,h$ of order $3$ and $2$ respectively,
acting on $\T$, $\R^2\cong \mathbb{C}$ and $\{1,2,3\}$
as follows:
\begin{eqnarray*}
&&  t\in \T \implies  g(t) = t, \ h(t) = t  \\
&&  v\in \mathbb{C} \implies
g(v) =  e^{\frac{2\pi}{3} J} v, \ h(v) = \bar v  \\
&& g(1) = 2,\ g(2) = 3,\ g(3) = 1,\ h(1) = 2,\ h(2) = 1,\ h(3) = 3.
\end{eqnarray*}
As a consequence, if $x(t)$ is an equivariant loop in $\Lambda^K$,
then, for each instant $t\in \T$, the configuration $x(t)$ 
is an (oriented) equilateral triangle  symmetric with respect 
to the real line. In other words, the configuration space 
has dimension $1$, with a subspace of dimension $0$
consisting of triple collisions.
Notice that $K$ acts trivially on the time circle $\T$.
Now extend $K$ by an (central) element $r$ of order $2$,
which acts on $\T$ by $r(t) = t+\pi$,
on $\R^2$ by $r(v) = -v$ and trivially on $\{1,2,3\}$.
It is easy to see that for the group $G = K \times C_2$ all
loops in $\Lambda^G$ have collisions at some time $t\in \T$.
\end{example}
Consider the normal subgroup $\ker \tau \triangleleft G$ 
(which will be termed 
\emph{core} of $G$ -- if it is trivial the group
will be said \emph{with trivial core})
\label{defi:core}\label{defi:trivialcore}
and the group $\bar G = G/\ker \tau$ acting on $\T$.
\begin{definition}
\label{action-types}
If the group $\bar G$ acts trivially on the orientation of $\T$ then
$\bar G$ is cyclic and we say that the action type of $G$ on $\Lambda$ is of
\emph{cyclic type}. If the group $\bar G$ consists of a single reflection on
$\T$, then we say that the action type of $G$ on $\Lambda$ is of 
\emph{brake type}. 
Otherwise we say that the action of $G$ on $\Lambda$ is of 
\emph{dihedral type}. 
\end{definition}

\begin{definition}
\label{T-isotropy}
The \emph{$\T$-isotropy subgroups of $G$} 
are the isotropy subgroups of the action of $G$ on the time
circle $\T$ (induced by $\tau$).
\end{definition}

\begin{definition}
\label{found_dom}
Let $\I\subset\T$ be the closure of a 
\emph{fundamental domain} 
for the action of $\bar G$ on $\T$ defined as follows.  When
the action type is cyclic then $\I$ is a closed interval connecting
the time $t=t_0$ with its image $zt$ with a cyclic generator $z$ in $\bar G$; in
this case $\I$ can be chosen among infinitely many intervals.
Otherwise, $\I=[t_0,t_1] \subset \T$ where $t_0,t_1$ are
distinct and have non minimal $\T$-isotropy subgroups in $G$;
moreover no other point in $(t_0,t_1)$ has non-minimal isotropy. There are
$|\bar G|$ such intervals.
\end{definition}

Let $H_0$ and $H_1$ be the isotropy subgroups of $t_0,t_1$ respectively. If
the action type is brake then $H_0=H_1=G$; while if it is dihedral, then $H_0$
and $H_1$ are distinct and proper subgroups of $G$. Further more the length of
$\I$ is always $\frac{T}{|\bar G|}$, since 
$\T =
\bigcup_{\bar g \in \bar G}\bar g\I$ 
and the interiors of the terms in the
sum are disjoint.

Let now consider the restriction of the Lagrangian action to the
$G$-equivariant loops
\begin{equation}
\label{G-A}
\action_\omega^G : \Lambda^G \rightarrow {\mathbb R}\cup \infty.
\end{equation}
\begin{equation*}
\action_\omega^G (x)=\int_{\T}L_\omega(x(t),\dot{x}(t))dt 
= |\bar G| \int_{\I}L_\omega(x(t),\dot{x}(t))dt.
\end{equation*}

The action functional~\ref{action:functional} is termed $G$-invariant when
$\action_\omega(x)=\action_\omega(gx)$; then the following result holds:
\begin{prop}[Palais Principle of symmetric criticality]
\label{PPSC}
If $\cal{A}_\omega$ is $G$-equivariant then a collisionless 
critical point of 
$\action_\omega^G$ 
is a critical point of $\action_\omega$.
\end{prop}   

\subsection{Coercivity and generalized solutions}
The action functional $\action_\omega^G: {\Lambda}^G \rightarrow \R \cup
\infty$ is termed
\label{ap:coercive}
\emph{coercive} in ${\Lambda}^G$ if $\action_\omega^G(x)$
diverges to infinity as the $H^1$-norm of $x$ goes to infinity in
${\Lambda}^G$. This 
is a fundamental property, for it guarantees the existence of
minimizers on the set ${\Lambda}^G$ for the functional $\action_\omega^G$.
Remark that this definition is equivalent to
\begin{lemma}
\label{existence} 
The action functional $\action_\omega^G(x)$ is coercive on ${\Lambda}^G$ if
and only if 
\[
\inf_{x \in \Lambda^G,  \left\| x \right\|_{L^2}=1} K_\omega >0.
\]
\end{lemma}
One  can easily deduce the following  corollary.
\begin{corollary}
\label{coro:coercive}
When $\omega \notin \ze$ (and the period $T=2\pi$),
for every finite group $G$ 
acting on the loop space
 the action functional 
$\action_\omega^G$ is coercive on ${\Lambda}^G$.
\end{corollary}
From now on we take the period $T$ to be $2\pi$,
in order keep the statements simple, unless otherwise explicitely
stated.

The next proposition immediately follows from 
Proposition  (4.1) in \cite{FT2003}.
\begin{prop}
\label{prop:coercive}
When $\omega= 0$, then  $\action_\omega^G = \action_0^G$ 
is coercive on ${\Lambda}^G$ if and only if $\mathcal{X}^G=0$. 
\end{prop}

When the functional $\action_\omega^G$ is coercive 
on the subspace
${\Lambda}^G$ and its minima are collisionless, the latter are
classical solutions of the Newton
equations~\ref{newton:eq}. To deal with the
possible occurrence of collisions 
some authors introduced the notion of \emph{generalized
solutions}. A generalized solution (of~\ref{newton:eq}) 
is a $H^1$ path $x(t)$
defined on an interval $(T_0,T_1)$ such that in $(T_0,T_1) \smallsetminus
x^{-1}\Delta$ is a classical 
$C^2$-solution, the Lagrangian action is
finite on $(T_0,T_1)$ and  partial energies are 
preserved after collisions:
\[
\forall t \in [t_0,t_1] \subset (T_0,T_1), 
\forall {\bf k}\subset{\bf n} : x(t) \notin \Delta_{{\bf k},{\bf k}'} 
\Longrightarrow E_{\bf k} \in H^1 ((t_0,t_1),\X),
\]
where $E_{\bf k}$  is the partial energy of the cluster ${\bf k}$. 

We can then conclude the study on the existence of minimizers of the
restricted action functional $\action_\omega^G$ with the following
proposition ((4.12) in \cite{FT2003})

\begin{propo} 
\label{prop:2.27}
When $\action_\omega^G$ is $G$-invariant and coercive on 
$\Lambda^G \subset \Lambda$, then there exists at least a 
minimum of the Lagrangian action $\action_\omega^G$ which yields a
generalized 
solution of~\ref{newton:eq} in ${\Lambda}^G$.
\end{propo}

Sometimes we shall refer the property of being coercive
directly to the symmetry group $G$.
\begin{defi}
\label{defi:Gcoervice}
A symmetry group $G$ is said to be \emph{coercive} if
the action functional $\action_0^G$, obtained with $\omega=0$,
is coercive 
in $\Lambda^G$.
\end{defi}
\section{Symmetry groups of type R and rotating frames}
\label{sec:typeR}

The determinant of the homomorphism $\rho\from G \to O(2)$
is defined by composition with $\det \from O(2) \to \{+1,-1\}$
(the same for $\tau$).
\begin{defi}
\label{def:typeR}
A symmetry group $G$ is said \emph{of type R} if the determinant
homomorphisms $\det(\rho), \det(\tau) \from G \to \{+1,-1\}$
coincide.
\end{defi}

\begin{lemma}
For every $\omega \in \mathbb{R}$, the Lagrangian action functional 
is invariant with respect
to the action of symmetry groups of type R.
\end{lemma}
\begin{proof}
The proof is straightforward.
\end{proof}

Now, consider the following change of coordinates
in $H^1(\R,\X)$ (\emph{rotating frame} with angular velocity $\omega$):
($\forall i \in \n$) $x_i(t) = e^{i\omega t} q_i(t)$.
A path $x(t)$ in $\Lambda$ is $G$-equivariant if an only if 
$x(gt) = gx(t)$ for every $g\in G$ and every $t\in \T$.
That means that for every $g\in G$ 
there exists a real number $\delta = \delta(g)$ 
such that  for 
the corresponding path $q(t)$ in the rotating frame
\begin{equation*}
e^{i\omega(t + \delta(g))} q(gt) = g e^{i\omega t} q(t)
\end{equation*}
for every $t \in \R$ if $g$ preservers the orientation in $\T$,
while if $g$ reverses the orientation 
\begin{equation*}
e^{i\omega(-t + \delta(g))} q(gt) = g e^{i\omega t} q(t)
\end{equation*}

If we assume $G$ to be of type R, 
then $g$ reverses the time orientation 
 if and only if it reverses the space orientation,
and therefore
in both cases
\begin{equation}
q(gt) = e^{-i\omega \delta(g)} g q(t).
\end{equation}
In order to let a single change of coordinates  
be enough to reduce as much as possible the group of symmetries,
one choses  $g$ as an element  in 
$\ker\det\tau \smallsetminus \ker\tau$ rotating
the time circle with minimal (non-zero) angle. Being a generator
of the (cyclic) quotient ${(\ker\det\tau)}/_{\ker\tau}$,
such an element and all its powers will act trivially on the plane,
once a suitable angular velocity $\omega$ is chosen. One can hence
prove the following useful lemma:

\begin{lemma}\label{lem:typeR}
Up to a suitable change of coordinates (in a rotating frame)
a symmetry group of type R has the following property: every 
symmetry $g$ acting as a non-trivial shift on the time line
acts as  the identity
on the euclidean space.
In other words, for every 
element $g\in \ker\det\tau \smallsetminus \ker\tau$ 
(equivalently,
every element $g\in G$ such that $\tau(g)$ 
is a  non-trivial
translation on the time line),
the image $\rho(g)$   is trivial:
$(\ker\det \tau \smallsetminus \ker\tau)\subset \ker\rho$.
\end{lemma}

\begin{figure}[ht]
\centering
\subfigure[Rotating frame]{%
\includegraphics[width=0.3\textwidth]{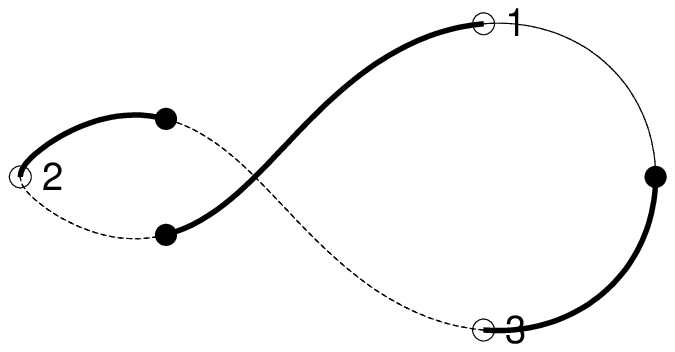}
\label{fig:notto_rf}
}
\hspace{0.15\textwidth}
\subfigure[Inertial frame]{%
\includegraphics[width=0.3\textwidth]{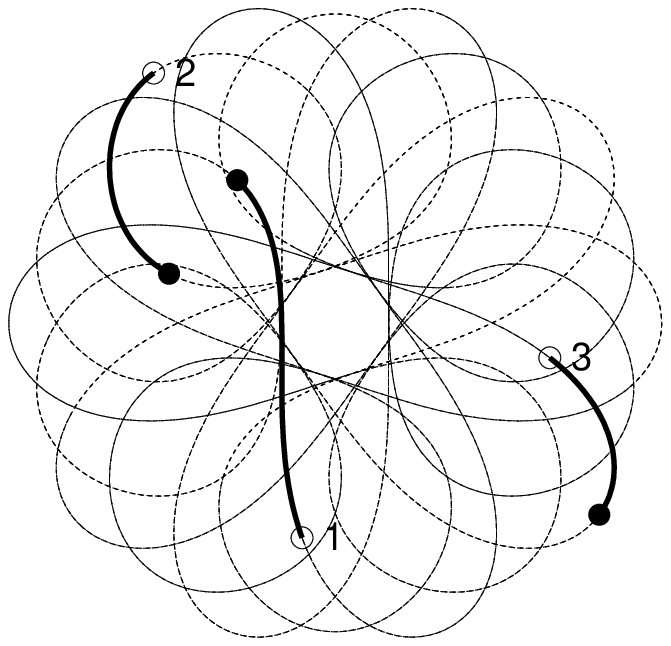}
\label{fig:notto}
}
\caption{A choreography symmetric with respect to a group of 
order $60$}
\end{figure}

Such a change of coordinates in rotating frames is useful for
simplifying orbits with many symmetries. For example, consider
the orbits in figure~\ref{fig:notto_rf} and \ref{fig:notto}
(these orbits  are equivariant minimizers 
for the group described in Remark~\ref{remark:notto}).
In the inertial (non-rotating) frame the orbit enjoys a symmetry
group of order $60$, with many rotations; the synchronization of the  rotating
frame with the symmetry rotations allows the symmetry group to be reduced 
to a much simpler form,
in rotating coordinates (practically, by canceling out the rotation
part in the symmetry group). This is just part of the wider
notion of \emph{redundant} symmetries. 

\begin{defi}
\label{defi:redundant}
An element $g$ of a 
symmetry group $G$ is 
termed \emph{redundant} if $g \neq 1$, 
$\tau(g) \neq 1$, $\det(\tau(g)) = 1$ 
and $\sigma(g) = \rho(g)  = 1$.  A group with redundant elements 
is termed a \emph{redundant} group.
\end{defi}

Hence a redundant symmetry is an element  $g$ 
which not only is 
an element $g\in G$ such that $\tau(g)$ is a non-trivial
translation on the time line and 
$\rho(g)$  is the identity, but also such that 
$\sigma(g)$ is the identical permutation.  Clearly redundant elements 
do not give any true constraint on the loops, and they can be modded out.
Given a symmetry group we will implicitely assume that the maximal redundant
subgroup (i.e. the subgroup consisting of all redundant symmetries) 
has been modded out, so that, unless otherwise stated, 
the group has no redundant elements. On the other hand, 
a change in rotating coordinates given as seen above can be 
used to create a redundant part that once mod out gives a smaller symmetry 
group. Summarizing, 
there are two ways of considering a group $G$ equivalent
to another group $G'$: either $G'$ is equal to $G/H$, where $H$ 
is a subgroup of redundant symmetries, or 
$G'$ is equal to $G$ after a change in the rotating angular 
velocity (for groups of type R). One might want to define an equivalence
relation in symmetry groups, and to consider always the \emph{reduced}
minimal group, that is the group with the minimum number of elements.
Since  sometimes it is easier to work with non-minimal groups,
we leave skip this notion of equivalence and will treat the groups
on a case-by-case fashion.

\begin{remark}\label{rem:typeR}
A symmetry group of type R
such that  all its elements
acting as a non-trivial shifts on the time line
act as  the identity
on the euclidean space (that is, a group
fulfilling the conclusion of~\ref{lem:typeR})
of course can be redundant, in a frame rotating with suitable
angular velocity $\omega$.
To see how $\omega$ changes when
we eliminate the redundant symmetries of the group
(and hence the period has to be rescaled) 
we proceed in the following way.
Let $g$ be an element in $\ker(\det(\tau)) \subset G$ which 
rotates the time circle $\T$ of a minimal angle (i.e.\
the elements inducing 
one of the two standard cyclic generators in $\ker(\det(\tau))/\ker\tau$).
Let $c$ denote its order.
If the group $G$ is of type R, then
$\rho(g)$ is a rotation of angle $2\pi b/c$, for some integer $b$ with $|b|<c$.
According to its definition, $\delta(g) = 2\pi/c$, 
and hence $\omega$ will be chosen in a way to satisfy the equation
\(\frac{\omega 2\pi}{c} = \frac{2\pi b}{c} \mod 2\pi\),
that is
\(\omega = b \mod c\).
On the other hand the true period of the paths in the rotating frame will be
a suitable multiple 
of 
$2\pi/c$ instead of $2\pi$ (depending on the period $s$
of the permutation $\sigma(g)$), hence, after making the action not redundant,
the corresponding $\omega$ will be simply $\omega' = \frac{s\omega}{c}$,
and thus chosen so that 
\begin{equation}\label{eq:omega}
\frac{\omega'}{s} = \frac{b}{c} \mod 1.
\end{equation}
Therefore a symmetry group of type $R$, not redundant and with 
$\omega'\neq 0$, in this method
yields back a symmetry group with $\omega =0$
provided that equation~\ref{eq:omega} can be solved for $b$, $c$
integers. Furthermore, the symmetry group found is unique
up to redundancies.
\end{remark}

For an equivariant path $x(t)$ 
let  $J(t)$ denote its angular momentum
\begin{equation}
\label{eq:angmom}
J(t) = 
\sum_{i\in \nn} m_i x_i \times \dot x_i
\end{equation}
(since we are dealing with planar paths, 
$J(t)$ belongs to a copy of $\R$ orthogonal to the plane $E$).

\begin{lemma}
\label{lemma:J}
For every equivariant $x(t) \in \Lambda^{G}$,
\[
J(gt)  = \det(\rho(g)) \det(\tau(g)) J(t) \in \R.
\]
\end{lemma}
\begin{proof}
We have the chain of equalities:
\begin{eqnarray*}
{J(gt)} & =  &
\sum_{i\in \nn} m_i x_i(gt) \times  \dot x_i(gt) 
\\
&=& 
\sum_{i\in \nn} m_i \left[ 
(g x_{g^{-1}i} (t)) \times  
( \det(\tau(g)) g \dot x_{g^{-1}i} (t)  )
\right] \\
&=& \det(\tau(g)) \sum_{i\in \nn} m_{g^{-1}i} \left[ 
(g x_{g^{-1}i} (t)) \times  
( g \dot x_{g^{-1}i} (t)  )
\right] \\
&= &
\det(\tau(g)) \sum_{i\in \nn} m_{g^{-1}i}  \det(\rho(g))  \left[ 
x_{g^{-1}i} (t) \times  
\dot x_{g^{-1}i} (t)  
\right] \\
&=& 
\det(\tau(g)) \det(\rho(g))  J(t). 
\end{eqnarray*}
\end{proof}
\begin{lemma}
\label{lemma:nottypeR}
If the symmetry group $G$ is not of type R, then
every $G$-equivariant path $x(t)$ has zero angular momentum.
\end{lemma}
\begin{proof}
Since $J$ is constant, 
if there is $g\in G$ such that $\det(\rho(g)) \neq \det(\tau(g))$
then
\[
J = \det(\rho(g)) \det(\tau(g)) J = - J 
\]
which means that $J=0$ if 
$G$ is not of type R.
\end{proof}
\begin{definition}
\label{def:fullyunc}
A symmetry group $G$ is termed \emph{fully uncoercive}
if it is neither of type R nor coercive (actually, this implies
that it is not possible to find a suitable $\omega$ 
such that the action is coercive -- on the converse, if the action
is of type R is always possible to find 
an $\omega$ such that the action functional
is coercive).
\end{definition}
\section{Planar symmetry groups}
\label{sec:class}
In this section we will list the principal symmetry groups 
for the planar $3$-body problem.

\subsection{The trivial symmetry}
\label{sym:trivial}
Let $G$ be the trivial subgroup of order $1$.
It is clear that it is of type R, it has 
the rotating circle property. It yields a coercive
functional on $\Lambda ^ G = \Lambda$ only 
for $\omega \neq 0 \mod 1$ (see~\ref{coro:coercive} 
and~\ref{prop:coercive}),
that is, if and only if $\omega$ is not an integer.
If $\omega = \frac{1}{2} \mod 1$ then the minimizers are minimizers
for the anti-symmetric symmetry group (which is also termed 
\emph{Italian} symmetry by some authors) $x(t+\pi) = -x(t)$.
It is worth noting that the masses can be different.
\begin{lemma}
\label{lem:actiontrivial}
For every $\omega \notin {\mathbb Z}$ and every choice of masses the minimum
for the trivial symmetry occurs in the relative equilibrium motion associated 
to the Lagrange central configuration.
\end{lemma}
\begin{proof}
The idea of the proof
comes from the proof of a similar proposition in 
\cite{chendesol}.
Consider the action functional defined 
in~\ref{action:functional} and the energy and potential functionals
\[
{\cal K}_\omega = 
\frac{1}{2} \sum_{i=1}^3 m_i \int_0^{2\pi} |\dot x_i-J\omega x_i|^2, \quad
{\cal U} = \int_0^{2\pi} U(x).
\]
The following estimate holds
\[
{\cal U} \geq U_0 \int_0^{2\pi} \frac{1}{I^{\alpha/2}},
\]
where $I=\sum_i m_i|x_i|^2$ is the momentum of inertia of the system 
and $U_0$ is the minimal value of the normalized potential $\tilde U(x)=
U\left(\frac{x}{I^{1/2}}\right)$. 

The equality ${\cal U}=U_0 \int_0^{2\pi} \frac{1}{I^{\alpha/2}}$ is 
reached if
and only if at every instant $t$ 
the configuration is
Lagrangian. Consider the Fourier representation of the trajectories
$x_i$, $x_i(t)=\sum_{n \in {\mathbb Z}} c_{i,n} e^{Jnt}$, $c_{i,-n}=\bar
c_{i,n}$, $i=1,2,3$; 
the kinetic part on the action functional can be estimated
as 
\[
{\cal K}_\omega = \frac{1}{2} \sum_{i=1}^3 m_i 
\int_0^{2\pi} \sum_{n \in {\mathbb Z}} (n-\omega)^2 c_{i,n}^2
\geq \frac{1}{2} \min_{n \in {\mathbb Z}}(n-\omega)^2  
\int_0^{2\pi} I = \frac{c(\omega)}{2} \int_0^{2\pi} I
\]
where $c(\omega) = (k-\omega)^2$ and  $k$ is the 
integer closest to $\omega$. 
The equality is attained if and only if 
$x_i(t)=\vec a_{i}\sin(kt)+J \vec b_{i}\cos(kt) $, 
where $\vec a_{i,k},\vec b_{i,k} \in \CC$, $i=1,2,3$.
We can conclude that
\[
\action_\omega(x) \geq  
\int_0^{2\pi} \frac{c(\omega)}{2}I + \frac{U_0}{I^{\alpha/2}}
\geq 2\pi \left( \frac{c(\omega)}{2}I_{min} + 
\frac{U_0}{I_{min}^{\alpha/2}} \right)
\]
for every loop $x \in \Lambda^G = \Lambda$. 
The equality is reached if and only 
if all the following conditions are verified:
\begin{enumerate}
\item\label{en:i}
the momentum of inertia is constantly 
equal to $I_{min}=\left( \frac{2U_0}{c(w)} \right)^{\frac{2}{\alpha+2}}$;
\item\label{en:ii} at every instant 
the configuration is Lagrangian;
\item\label{en:iii} 
there exist $\vec a_{i},\vec b_{i} \in {\mathbb R}^2$, 
such that 
$x_i(t)=\vec a_{i}\sin(kt)+J \vec b_{i}\cos(kt)$, $i=1,2,3$, 
where $k$ is the integer closest to $\omega$.
\end{enumerate}
From~\ref{en:ii} 
we deduce the existence of a 
time-dependent complex-valued function $\lambda(t)$ such that 
$x_i(t)=\lambda(t)\xi_i$, where $(\xi_i)_i$ is a Lagrange central 
configuration. 
But the momentum of inertia is fixed and, up to
reflection there is  a unique central configuration
with $I=I_{min}$. The complex
function $\lambda$ is then a rotation 
fixing the center of mass of the system
and the trajectories $x_i(t)$ are circles.  
\end{proof}
\subsection{The line symmetry}
\label{sym:line}
Another case of symmetry group of type R with
arbitrary masses is the line symmetry:
the group is a group of order $2$ acting
by a reflection on the time circle $\T$,
by a reflection on the plane $E$,
and trivially on the set of indexes $\nn$.
That means, at time $0$ and $\pi$ the masses are
collinear, on a fixed line $l\subset E$.
By~\ref{coro:coercive} and~\ref{prop:coercive}
it is coercive only when $\omega \not\in \ze$.
In this case the Lagrangian solution cannot be a minimum.

\begin{figure}[ht]
\begin{center}
\psfrag{action}{$\action_\omega$}
\psfrag{omega}{$\omega$}
\includegraphics[width=0.7\textwidth]{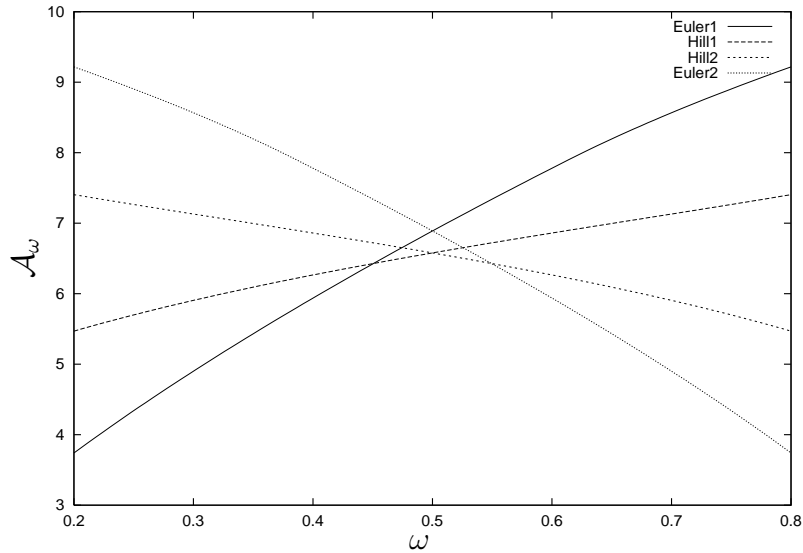}
\caption{Action levels for the line symmetry}
\label{fig:line}
\end{center}
\end{figure}

Consider the homographic motions from the Euler configuration 
(with the third body in the center)
\begin{equation}
\label{EO}
x_1(t)=R e^{Jkt},  x_2(t)=-R e^{Jkt},  x_3(t)=0 \quad 
(\mbox{for } t \in [0,2\pi], R>0,  k \in \ze).
\end{equation}
These trajectories are the {Euler orbits}.
In a rotating frame with angular velocity $\omega$ (when the masses
are all equal to $1$) 
the action functional can be computed as
\begin{equation*}
\action_E (R,\omega,\alpha) 
= R^2 (k-\omega)^2 + \frac{2}{R^\alpha} +\frac{1}{(2R)^\alpha}.
\end{equation*}
Thus, for $\alpha=1$, 
the minimal value of $\action_E$ is attained
on $R_\omega= \sqrt[3]{\frac{5}{4(k-\omega)^2}}$ and
\begin{equation}
\label{minactionEO}
\min_{R>0} \action_E (R,\omega,1) =  2\pi \frac{3}{2}\frac{\sqrt[3]{25(k-\omega)^2}}{\sqrt[3]{2}}.
\end{equation}
In particular when $\omega = \frac{1}{2}$, 
\begin{equation}
\label{eq:mineuler}
\min_{R>0} \action_E (R,\frac{1}{2},1)= 
\action_E (\sqrt[3]{5},\frac{1}{2},1) 
= 2\pi\frac{3}{4}\sqrt[3]{25}.
\end{equation}
Actually, there are two minima, one for $k=0$
and one for $k=1$ (the action levels of these
two Euler orbits can be found in figure~\ref{fig:line},
under the names {\tt Euler1} and {\tt  Euler2}).

Our aim now is 
to prove the existence of a motion which is not homographic,
whose action level is less than 
the minimum on the Euler orbits, 
for $\alpha=1$ and  $m_i=1$ for  $i=1,2,3$. 
Consider the following family of test paths
\begin{equation}
\label{HO}
x_1(t)=d+Re^{ikt}, x_2(t)=d-Re^{ikt}, x_3(t)=-2d,
\end{equation}
for $t \in [0,2\pi]$, $0<R<3d$,  $k \in \ze$.
The center of mass lies in the origin. 

The action functional is
\begin{eqnarray}
\label{eq:testpath}
\frac{\action_H (R,d,\omega)}{2\pi} &=& 
\frac{1}{2\pi}\int_0^{2\pi} 
\left( |\dot x_1 -J\omega x_1|^2 + 
\frac{4\omega^2 d^2}{2} + 
\frac{1}{2R} +\frac{2}{\left|3d+Re^{Jkt}\right|} \right) \\
\notag &=& 
3\omega^2 d^2 + R^2(k-\omega)^2 + \frac{1}{2R}
+\frac{2}{3d} \frac{1}{2\pi} 
\int_0^{2\pi} \frac{dt}{\left|1+\frac{R}{3d}e^{Jkt}\right|}.
\end{eqnarray}
Thus when $\omega=\frac{1}{2}$,
by Lemma (8.7) and (8.13) in \cite{FT2003},
\begin{equation}
\label{hill_omega}
\frac{\action_H (R,d,\frac{1}{2})}{2\pi} \leq 
\frac{3}{4}d^2 + \frac{1}{4}R^2 + \frac{1}{2R}
+\frac{2}{3d} \left( 1-\frac{1}{2}\log(1-\frac{R^2}{9d^2}) \right).
\end{equation}
Now observe that 
by setting
$R=1$ and $d=\frac{4}{5}$ one obtains
\begin{eqnarray*}
\frac{\action_H (1,\frac{4}{5},\frac{1}{2})}{2\pi} &\leq& 
\frac{3}{4}\frac{4^2}{5^2} + \frac{1}{4} + \frac{1}{2}
+\frac{5}{6} \left( 1-\frac{1}{2}\log(1-\frac{5^2}{9\cdot 4^2}) \right) \\
&=&  \frac{619}{300} - \frac{5}{12} \log(\frac{119}{144}).
\end{eqnarray*}
Hence, by~\ref{eq:mineuler} the claim follows once it is shown that
\begin{equation*}
\frac{619}{300} + \frac{5}{12} \log(\frac{144}{119}) 
< \frac{3}{4}\sqrt[3]{25}.
\end{equation*}
But one has $\frac{144}{119} = 1 + \frac{25}{119}$ 
and hence
\begin{equation*}
\frac{619}{300} + \frac{5}{12} \log(\frac{144}{119})  <
\frac{619}{300} + \frac{5}{12} \frac{25}{119}  
=
\frac{38393}{17850}
\end{equation*}
which can be easily (and rigorously)
seen to be less than $\frac{3}{4}\sqrt[3]{25}$.

In figure~\ref{fig:line}, besides the action levels
of two  Euler orbits,
it is shown the graphs of the action functional of two minimizer
orbits (computed numerically): 
the curve there named {\tt Hill1} corresponds to 
a test path with $k=1$ and
the curve named {\tt Hill2} corresponds to the value
$k=0$.
By continuity of 
$\action_H$ with respect to the masses 
$m_i \in {\mathbb R}^*_+$, $i=1,2,3$, and to $\omega$, 
the following result holds.

\begin{lemma}
\label{le:LS}
There exists an open set containing $\frac{1}{2}+\ze \subset \R$
such that if the masses are approximately equal 
then the minimum of the action 
functional restricted to the equivariant orbits with respect to a 
line symmetry group is not homographic (i.e.\ not the relative 
equilibrium motion associated to the Euler configuration).
\end{lemma}

\begin{figure}[ht]\centering
\subfigure[Rotating frame]{%
\includegraphics[width=0.3\textwidth]{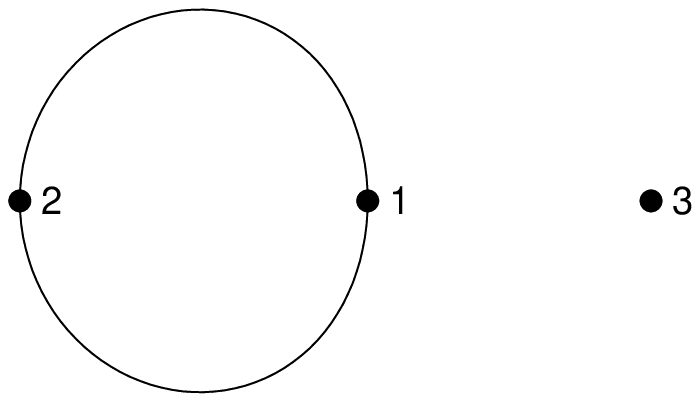}
\label{fig:sol_hill}
}
\hspace{0.15\textwidth} 
\subfigure[Inertial frame]{%
\includegraphics[width=0.3\textwidth]{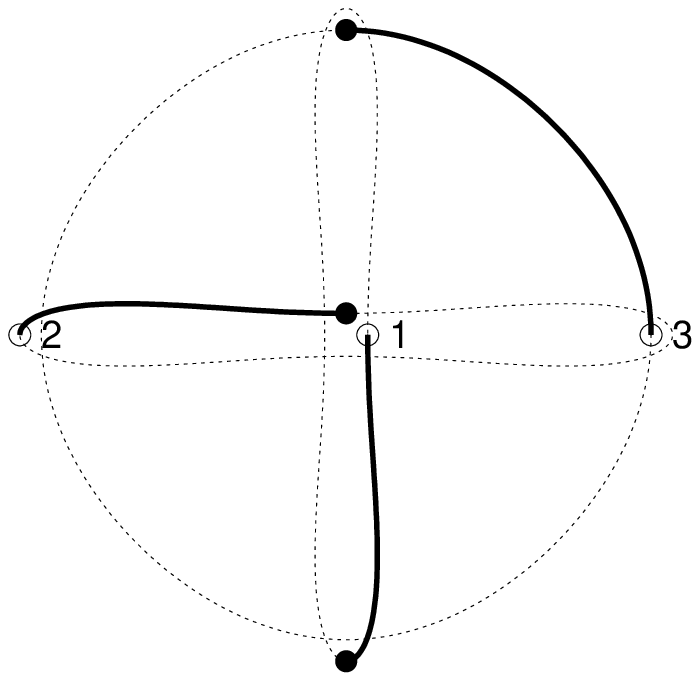}
\label{fig:nicesol}
}
\caption{Hill-like  orbit (equal masses and $\omega=1$)}
\end{figure}

A possible candidate for such minimizer is depicted in 
figure~\ref{fig:sol_hill} (in rotating coordinates)
and in figure~\ref{fig:nicesol} (inertial coordinates). 
It is a well-known orbit, which was found
numerically for instance in \cite{henon,moore}.
In particular, Chen in \cite{chen2003_II}
propose a proof for this orbit in which numerical methods
play a very minor role and are carefully justified.
The reader is referred to \cite{or} (which is also
available on {\tt www.arxiv.org} or from the web-site
of the MPI) and to \cite{chen2003_II,chen2003_I} 
for other orbits and symmetries of this type. 
Numerical action level estimates  of retrograde test paths
can be found also
in 
section 4 of \cite{chen2003_I} and
section 7 of \cite{or}.
 
\subsection{The $2$-$1$-choreography symmetry}
\label{sym:21}
Consider the group of order $2$ acting as follows:
$\rho(g) = 1$, $\tau(g) = -1$ (that is, 
the translation of half-period) and $\sigma(g) = (1,2)$
(that is, $\sigma(g)(1) = 2$, $\sigma(g)(2) = 1$, and $\sigma(g)(3) = 3$.
That is, it is a half-period  choreography for the bodies $1$ and $2$.
It is of type R,  and coercive for a suitable choice of $\omega$.
By~\ref{prop:coercive} and~\ref{coro:coercive}, 
the choice is $\omega \notin \ze $. The first two
masses have to be equal $m_1=m_2$.
Consider as above the case of equal masses.
Euler orbits (defined  in~\ref{EO}) and 
the test paths (defined in~\ref{HO})  are $G$-equivariant 
if and only if $k=1 \mod 2$ (this, incidentally,
implies that $\action_\omega^G$ is coercive if and only
if $\omega\not\in \ze$). 
Roughly speaking, this implies that the branch 
termed {\tt Euler1} in figure~\ref{fig:line} (which
is a Euler path with $k=0$) does not appear in a corresponding
graph for $2$-$1$-choreographies, while both {\tt Hill1}
and {\tt Hill2} do.
We are going to prove that actually this is the case,
and that therefore under such a constraint 
the minimizer is not homographic for a rather bigger
interval of angular speeds $\omega$.
\begin{lemma}
\label{le:21C}
There is an open set $A$  containing $[-\frac{1}{2},\frac{1}{2}]$ such that 
if $\omega \in A+2\ze$
then minimizers under  the $2$-$1$-choreography symmetry group
are not homographic.
(i.e.\ not the relative 
equilibrium motion associated to the Euler configuration).
\end{lemma}

\begin{figure}[ht]
\begin{center}
\psfrag{action}{$\action_\omega$}
\psfrag{omega}{$\omega$}
\includegraphics[width=0.7\textwidth]{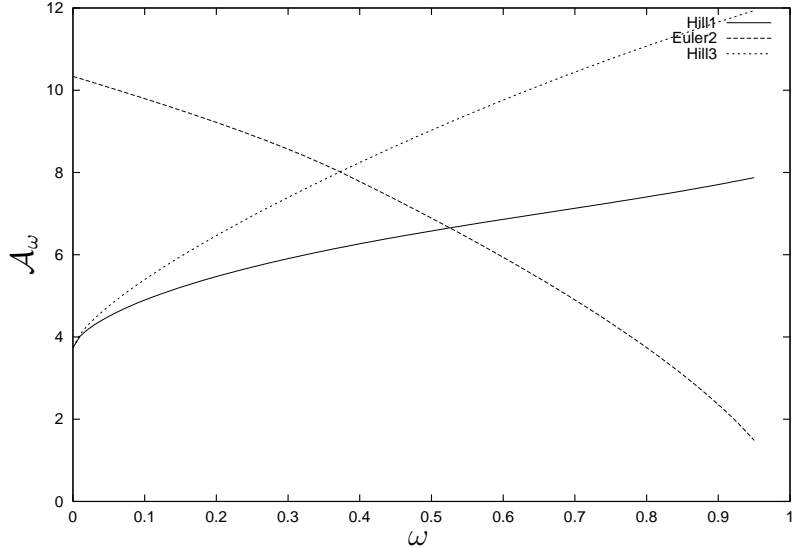}
\caption{Action levels for the $2$-$1$-choreography symmetry}
\label{fig:choreo_21}
\end{center}
\end{figure}

\begin{proof}
Consider $\omega \in [0,\frac{1}{2}]$; 
the action $\action_E$ of an Euler minimizer is given 
by~\ref{minactionEO} with $k=1$, since for 
$k=0$ the orbit is not equivariant.
This is a decreasing function of $\omega \in [0,\frac{1}{2}]$.
On the other hand, the action functional $\action_H$
evaluated at  the test path~\ref{HO}, as in~\ref{eq:testpath} (with $k=1$)
is a (convex) function of $\omega$.
It has already been shown in~\ref{le:LS}
that, for $\omega=\frac{1}{2}$, $\action_H<\action_E$.
Deriving the difference
$\action_H - \action_E$ after substituting $d=\dfrac{4}{5}$ 
and $R=1$ as above 
one obtains a function 
\[
\dfrac{d}{d\omega}(\action_H - \action_E) = 
2\pi \left( 
\dfrac{146}{25} \omega - 2
+ (\frac{25}{2})^{1/3} (1-\omega)^{-1/3}\right),
\]
which can be easily seen to 
be increasing in $[0,\frac{1}{2}]$ and positive 
in $0$. Thus the difference $\action_H - \action_E$ is 
increasing in $[0,\frac{1}{2}]$, and hence everywhere negative 
since it is so in $\dfrac{1}{2}$.
The graphs of the action levels are shown in figure~\ref{fig:choreo_21},
under the labels {\tt Hill1} and {\tt Euler2}.
\end{proof}

\begin{remark} 
We give a few words on 
figure~\ref{fig:choreo_21},
which represents the action levels (computed numerically).
There is a third branch of 
action levels, labeled {\tt Hill3}, which corresponds 
to a non-retrograde Hill-like orbit. 
In fact, for $\omega=0$, it is easy to see that the 
two action levels correspond to a Kepler $2$-body problem (where,
as $\omega$ tends to $0$, the third particle escapes
to infinity -- the functional
is not coercive) and a rotating Euler solution with 
the third particle in $0$.  More precisely, for $\omega=0$
there are two $2$-body 
Kepler solutions (one rotating clockwise with $k=-1$  and 
the other counter-clockwise with $k=1$) and two Euler solutions
(one clockwise  with $k=-1$ and the other counter-clockwise with $k=1$) 
as limit cases.
As $\omega$ increases (thus considering a frame rotating counter-clockwise),
the action of the Euler solution will decrease if $k=1$
and will increase if $k=-1$  (the graph of the latter Euler
solution does not appear in figure~\ref{fig:choreo_21}).
On the other hand, both the Hill-like solutions will
have an increased action, since there will 
appear the interaction of the third body, 
but the action of the orbit with $k=1$ turns out to be less than that
with $k=-1$, as one would expect. 
A possible trajectory for the minimizer (with $\omega=\frac{2}{5}$ and 
equal masses) in the inertial frame is shown in figure~\ref{fig:five}.
\begin{figure}[ht]
\begin{center}
\includegraphics[width=0.3\textwidth]{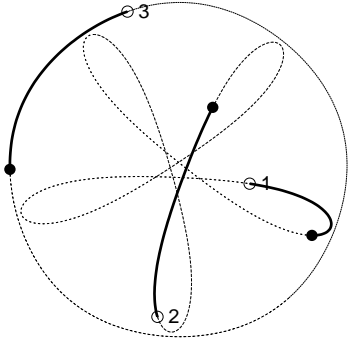}
\caption{A $2$-$1$-choreography with $\omega=\frac{2}{5}$}
\label{fig:five}
\end{center}
\end{figure}
Equivalently, this orbit can be found as a minimizer of $\action_0^G$,
where the redundant symmetry group $G$ is cyclic of order $10$ 
(hence the action is of cyclic type)  and 
generated by an element $g$ which 
acts as follows: $\tau(g)$ is a cyclic rotation
of $\pi/5$, $\rho(g)$ is a rotation of angle $2\pi/5$ and
$\sigma(g)$ is the permutation $(1,2)$. 
More generally, choosing $\rho(g)$ to be a rotation
of angle $\theta 2\pi$, with $\theta$ rational, yields
a non-homographic 
periodic symmetric orbit provided that $\theta\in(0,\dfrac{1}{4}]$.
\end{remark}

\subsection{The isosceles symmetry}
\label{sym:isosceles}
The isosceles symmetry can be obtained as follows:
the group is of order $2$, 
generated by $h$; 
$\tau(h)$ is a reflection in the time circle $\T$,
$\rho(h)$ is a reflection along a line $l$ in $E$,
and $\sigma(h) = (1,2)$ as above.
The constraint is therefore that at time $0$ and $\pi$
the $3$-body configuration is an isosceles triangle
with one vertex on $l$ (the third).
By~\ref{coro:coercive} and~\ref{prop:coercive}
it is coercive if and only if $\omega \not\in \ze$.
\begin{lemma}
For every $\omega \not\in \ze$ and every choice of masses (compatible
with the isosceles constraint) 
the minimum
for the isosceles symmetry occurs 
in the relative equilibrium motion associated 
to the Lagrange configuration.
\end{lemma}
\begin{proof}
It follows directly
from~\ref{lem:actiontrivial}.
\end{proof}

\subsection{The Euler--Hill symmetry}
\label{sym:hill}
Now consider the symmetry group  with a cyclic
generator $r$ of order $2$ (i.e.\ $\tau(r)=-1$) and a time 
reflection $h$ (i.e.\ $\tau(h)$ is a reflection of $\T$) given by
$\rho(r) = 1$, $\sigma(r) = (1,2)$,
$\rho(h)$ is a reflection and $\sigma(h)=()$.
It contains the $2$-$1$-choreography (as the subgroup $\ker\det(\tau)$), 
the isosceles symmetry (as the isotropy of $\pi/2\in \T$) 
and the line symmetry (as the isotropy of $0\in \T$)
as subgroups. 

Since Euler orbits (defined in~\ref{EO}) and the Hill-like test paths
(defined in~\ref{HO}) are 
equivariant if $k=1$, one can then draw the following immediate
consequence 
(see~\ref{le:LS} and~\ref{le:21C}).

\begin{lemma}
The minimum of the Euler--Hill symmetry is not homographic
for all values of $\omega$ in an open set containing 
$[-\frac{1}{2},\frac{1}{2}] + 2\ze$,
provided that the masses are approximately equal and 
two of them equal.
\end{lemma}

\subsection{The choreography symmetry}

\label{sym:choreography}
The choreography symmetry is given by the group $C_3$ of order $3$
acting trivially on the plane $E$, by a rotation 
of order $3$ in the time circle $\T$ and 
by the cyclic permutation $(1,2,3)$ in $\nn$.
The proof of~\ref{lem:actiontrivial} can be modified in a few
points in order to give the following result, which has been 
proved in   \cite{BT2004} -- Remark 4. 
It is easy to see that  the action 
functional is coercive if and only
if $\omega\neq \pm 1/3 \mod 1$.
\begin{lemma}
\label{lem:actionchoreo}
For every $\omega \neq \pm 1/3 \mod 1$ 
the minimal choreography of the three--body problem 
is a rotating Lagrange configuration.
\end{lemma}

\subsection{The Lagrange symmetry}
\label{sym:lagrange}
The Lagrange symmetry group  is the extension of the choreography 
symmetry group  by the isosceles symmetry group. Thus, 
it is a dihedral group of order $6$ and  the action is of type R.
Hence, the relative equilibrium motions associated 
to the Lagrange configuration are admissible motions for this symmetry
and the following result follows from 
both~\ref{lem:actiontrivial} and~\ref{lem:actionchoreo}.
As for the choreography symmetry, 
the action functional is coercive if and only
if 
$\omega\neq \pm 1/3 \mod 1$.
\begin{lemma}
For every $\omega \neq \pm 1/3 \mod 1$ 
and every choice of masses the minimum
for the Lagrange symmetry occurs in the relative equilibrium motion associated 
to the Lagrange configuration.
\end{lemma}

\subsection{The Chenciner--Montgomery symmetry group and the eights}
\label{sym:eights}

In this section we will describe 
three symmetry groups (up to change of coordinates)
that might 
yield the Chenciner--Montgomery \emph{figure eight} orbit as a minimizer:
they are the only symmetry groups not of type R for $n\leq 3$
bodies in the plane.

First, we consider the only symmetry group of cyclic action type,
not of type R and coercive (under these hypotheses the group 
is automatically not bound to collisions and transitive, 
as we will see later).
The group $C_6$ (the \emph{cyclic eight}
has order $6$, acts cyclically on $\T$ (i.e.\
by a rotation of angle $\pi / 3$),
by a reflection in the plane $E$, and by
the cyclic permutation $(1,2,3)$ in the index set.
The second group, which we denote by $D_{12}$,
is the group of order $12$ obtained by
extending $C_6$ with the element $h$ defined 
as follows: $\tau(h)$ is a reflection in $\T$,
$\rho(h)$ is the antipodal map in $E$ (thus, the rotation
of angle $\pi$), and $\sigma(h)$ is the permutation
$(1,2)$. This is the symmetry group used in \cite{monchen}.

The third group is the subgroup of $D_{12}$
generated by $h$ and the subgroup $C_3$ of order $3$ of $C_6 \subset D_{12}$.
We denote this group $D_6$ (since it is a dihedral group of order $6$).
The symmetry groups $D_{12}$  and $D_6$ are of dihedral type.
The choreography group $C_3$ is a subgroup of all the three 
groups, thus they are coercive.
Since they are not of type R, the minimum is not a homographic solution
(since it has trivial $J$).

\begin{propo}
\label{prop:D6D12}
If $\alpha>2$, then 
the minimum in $\Lambda^{D_{12}}$ (the Chenciner--Montgomery eight) 
coincides with the minimum
in $\Lambda^{D_6}$:
Any local  minimizer symmetric with respect to the group 
$D_6$ is symmetric with respect to the group $D_{12}$.
\end{propo}
\begin{proof}
Let $H_0,H_1 \subset D_6$ be the subgroups 
generated respectively by the elements $h_0$ and $h_1$ given
by $\rho(h_0) = \rho(h_1) = -1$
and $\sigma(h_0) = (1,2)$, $\sigma(h_1) = (2,3)$.
If $x(t)$ is a $D_6$-equivariant,  then
$t\in \T^{H_0} \implies x(t) \in \X^{H_0}$
and $t\in \T^{H_1} \implies x(t) \in \X^{H_1}$.
Thus, without loss of generality we can assume that 
at time $t=t_0=0$ one has $x(t_0) \in \X^{H_0}$,
while at time $t_1= \frac{2\pi}{6}$ one has
$x(t_1) \in \X^{H_1}$.
Since $[t_0,t_1]$ is a fundamental domain for the action
of $D_6$  in $\T$, if $x(t)$ is a $D_6$-equivariant (local)
minimum of the action functional $\action_\omega$,
then it minimizes the action restricted to the interval
$[t_0,t_1]$ together with  the constraints
$x(t_0) \in \X^{H_0}$ and $x(t_1) \in \X^{H_1}$.
Now consider the function $f\from \X \to \R$
defined by
$f(x_1,x_2,x_3) = (x_1-x_3)x_2$.
If a configuration (without collisions) 
$(\xi_1,\xi_2,\xi_3)$ belongs
to $\X^{H_0}$, then 
$\xi_1 = -\xi_2$ and $\xi_3 = 0$,
therefore 
$f(\xi_1,\xi_2,\xi_3) = - \xi_1^2 < 0$. 
On the other hand, if $(\xi_1,\xi_2,\xi_3)$ belongs
to $\X^{H_1}$ then
$f(\xi_1,\xi_2,\xi_3) = \xi_2^2 >0$.
It follows that there has to be a time $t'\in (t_0,t_1)$
such that $f(x(t'))=0$, that is,
such that the triangle given by the configuration
$x(t')$ is isosceles. Let $l$ be the reflection in the plane
such that $x_1(t') = l x_3(t')$, $x_2(t') = lx_2(t')$.
Let $g$ be the element of $O(\T)\times O(d) \times \Sigma_3$
given by: 
$\rho(g) = l$, $\sigma(g) = (1,3)$ and $\tau(g)$ 
is the time reflection fixing the time $\frac{\pi}{6}$.
It is easy to see that $g \X^{H_0} = \X^{H_1}$, 
$g\X^{H_1} = \X^{H_0}$ and $gx(t') = x(t')$.
Now, the restriction of $x(t)$ to the interval 
$[t_0,t']$ is not necessarily a (local) minimizer in the class 
of paths defined in $[0,t']$ with constraints
$x(0) \in \X^{H_0}$ and $x(t') \in \X^g$,
but sure its action is bounded below by such a minimal value.
The same holds for the restriction of $x(t)$ to 
the interval $[t',t_1]$, 
with respect to paths
with constraints $x(t') \in \X^g$ and $x(t_1) \in \X^{H_1}$.
Now, if $\action_1$ denotes the corresponding local minimum
of the action on paths defined in $[0,1]$
with boundary constraints $x(0) \in \X^{H_0}$ 
and $x(1) \in \X^g$,
then the sum of the actions over $[t_0,t']$
and $[t',t_1]$ is 
equal to 
\[
\left( T_0^{\frac{2-\alpha}{2+\alpha}} +
 T_1^{\frac{2-\alpha}{2+\alpha}} \right)
\action_1,
\]
where $T_0 = t'$ and $T_1=t_1-t'$.
Such a sum is minimal for $T_0 = T_1$,
hence 
the action of $x(t)$ (the minimal 
path with constraints $x(t_0) \in \X^{H_0}$ and 
$x(t_1) \in \X^{H_1}$) in $[t_0,t_1]$
is greater than 
twice the action of the path 
minimizing with constraints
$x(t') \in \X^g$ and $x(t_0) \in \X^{H_0}$ 
where $t' = \frac{\pi}{6} = \frac{1}{2}(t_1 - t_0)$.
But this minimizer, say $y(t)$, can be extended to 
a path defined on $[t_0,t_1]$
by the action of $g$:
\[
y(t'+t) = g y(t'-t).
\]
Hence $x(t)$ has to be equal to the symmetric path $y(t)$
on $[t_0,t_1]$. Now, the group generated
by $D_6$ and $g$ is conjugated to $D_{12}$, 
and this implies that any $D_6$-symmetric local minimizer
has the full $D_{12}$ symmetry group.
\end{proof}
\begin{remark}
Numerical experiments suggest that also the minimum 
in $\Lambda^{D_6}$ might coincide with 
the minimum in $\Lambda^{C_6}$  (and therefore that the 
occurring eight solution is unique), and that the restriction
on the exponent is not necessary.
Furthermore, 
it is not difficult to generalize~\ref{lemma:J} to orbits
in $3$-dimensional space, and thus to find symmetry groups
for which all minimizers have zero angular momentum.
Since $3$-body orbits whose
angular momentum vanishes are planar, 
in this paper actually we provide a complete
list of symmetry groups (in the space as in the plane)
not of type R, and hence that even in space eight-type orbits
are at most two (one for the $C_6$ and one for $D_{12}$). 
It is interesting to find (if possible) a proof of the fact
that actually it is unique.
\end{remark}

\begin{remark}\label{remark:notto}
As it will be proved in the next section,
we have described all possible symmetry groups 
for the planar three-body problem (in the sense of 
section~\ref{sec:symmcon}).
It is interesting that,
even when the global minimum is proved to be homographic, 
such symmetry groups can be used to find (at the moment
only numerically, but after estimates on the action
levels of colliding trajectories rigorous results might follow)
symmetric periodic orbits which are just local minimizers,
such as the orbit of figure~\ref{fig:notto}.
This orbit is a choreography obtained as local minimum
for the Lagrange symmetry group with $\omega=\frac{3}{10}$,
which corresponds to a dihedral  symmetry 
group $G$ of order $60$ acting as follows: if $g$ 
denotes the generator of the maximal cyclic subgroup  of $G$
and 
$h$ one of the elements of order $2$ not in the center,
then $\rho(g)$ is a rotation of angle $\frac{\pi}{5}$,
$\sigma(g) = (1,2,3)$, 
$\rho(h)$ is a reflection along a line,
$\sigma(h) = (1,2)$,
$\tau(g)$ is a rotation of angle $\frac{\pi}{15}$
and $\tau(h)$ a time-reflection.
In fact, for $\omega=0$ the figure eight orbit
seems to be a local minimum for the Lagrange symmetry
group, and this minimum can be continuated as $\omega$ grows
to obtain a family of quasi-periodic orbits, which 
are choreographies in the rotating frame. It is easy to see 
that there is a dense subset $D$ of $(0,1)$ 
such that for $\omega\in D$ the corresponding 
orbit is not only a choreography in the rotating frame,
but also a choreography in the inertial frame.
\end{remark}

\section{The classification}
\label{sec:proof}

In this section we will prove Theorem~\ref{MT1}. It will be an immediate 
consequence of~\ref{prop5.2} and~\ref{nr:nontrivialcore}.
 
\begin{lemma}
\label{lemma:dec}
If $G$ is a symmetry group with trivial core 
for the $3$-body problem in the plane,
and the action is not redundant,
then (up to a rotating frame) $\rho(g^2) = 1$ for every $g\in G$.
\end{lemma}
\begin{proof}
Let $r\in G$ denote the generator of $\ker\det(\tau)$
and $h$ a time reflection in $G$ (if it exists).
It must be that $hrh^{-1}=r^{-1}$ and $h^2=1$. 
Thus the conclusion is equivalent to  $\rho(r^2)=1$.
If the action is of type R and not redundant, 
then by \ref{lem:typeR} one has $\rho(r)=1$.
So assume that the action is not of type R.
If the action type is cyclic 
then $\rho(r)$ is a reflection,  so that 
$\rho(r^2)=1$. 
If the action type is brake, then 
one has $\rho(r)=1$ so the same follows.
If the action type is dihedral and 
$\rho(r^2)\neq 1$, then $r\in\ker\det(\rho)$
(since elements not in $\ker\det(\rho)$ have order
at most $2$). Thus, being the action not of type R,
$\rho(h)$ is orientation-preserving and of order at most $2$, 
i.e.\ $\rho(h)$ is either trivial or the antipodal map 
and hence commutes with $\rho(r)$.
Thus, by the equation $hrh^{-1} = r^{-1}$
it follows that
$\rho(r) = \rho(r^{-1})$
and therefore that $\rho(r^2) = 1$ as claimed.
\end{proof}

\begin{propo}
\label{prop5.2}
For the planar $3$-body problem with trivial core
the symmetry groups  not bound to collisions, not redundant
and 
not fully uncoercive
are 
those described in section~\ref{sec:class} and summarized
in table~\ref{table:class}.
\end{propo}
\begin{proof}
We can begin with symmetry groups acting trivially on the index set.
First we consider groups of type R:
any group of cyclic action type and not redundant acts trivially
on $E$, so we have to assume $G$ to be 
of brake or dihedral type, if it is not trivial
(see~\ref{sym:trivial}). But for the same reason
it cannot be of dihedral type (since the cyclic part would be 
trivial), hence only brake.
Let $h$ denote the time reflection. Since $G$ is of type R,
then $\rho(h)$ has to be a reflection in $E$,
and hence the group is the symmetry group of the line symmetry (see~\ref{sym:line}).
Now, consider a group not of type R. We want to show that 
is always either redundant or bound to collisions.
Let $r$ and $h$ be its cyclic generator (i.e.\ the generator
of $\ker\det(\tau)$) and one of its time reflections (i.e.\
one of the elements not in $\ker\det(\tau)$).
If the action type is cyclic, 
then $\rho(r)$ is  a reflection, 
and $G$ is not coercive.
If the action type is brake, 
then $\rho(h)$ is either trivial (and therefore $G$ is not coercive)
or the antipodal map (and therefore bound to collisions).
On the other hand if the action type is dihedral, 
then $\rho(r)$ is a reflection along a line $l$ 
in $E$ (otherwise would be trivial
or redundant). As a consequence, 
if $\rho(h)$ is trivial or a reflection along the same
line $l$, then $G$ is not coercive.
The other possibility is that $\rho(h)$ is a reflection
along a line $l'$ orthogonal to $l$, but in this case 
$G$ would result to be bound to collisions.
Thus there are no symmetry groups with the required properties
and not of type R.

Now we come to 
symmetry groups with transitive decomposition $2+1$.
Let $r$ be the cyclic generator  of 
$\ker\det(\tau)$ and, if the action type is not
cyclic,  $h$ an element of $G\minus \ker\det(\tau)$ (a time
reflection). 
First consider groups of type R. If the action type
is cyclic, 
then, to avoid being redundant, the only choice
is $\rho(r)=1$ and  $\sigma(r)=(1,2)$, 
thus the $2$-$1$-choreography (see~\ref{sym:21}).
If the action type is brake, 
then $\rho(r)$ and $\sigma(r)$ are trivial and thus $\rho(h)$ is a reflection
(since the group is of type R),
while $\sigma(h) = (1,2)$ (up to a permutation of the indexes). 
This is the isosceles symmetry group (see~\ref{sym:isosceles}).
As a third case, assume that 
the action type is dihedral. 
Being not redundant, it must be that 
$\rho(r)=1$, so that $\sigma(r) = (1,2)$ 
(otherwise it would be brake);
as above $\rho(h)$
is a reflection while $\sigma(h) = () $ or $\sigma(h) = (1,2)$.
Both cases yield the same symmetry group (the two
groups are conjugated in $O(\T)\times O(2) \times \Sigma_3$), 
which is the Hill symmetry group (see~\ref{sym:hill}). 
It can be generated as the union of any two  
of the three listed symmetry groups of order $2$.

Then, we deal with groups not of type R.
An enumeration of all possible cases shows that: 
\begin{nr}\label{nr:temp1}
There are no symmetry groups not of type R 
which are coercive, with transitive 
decomposition of type $2+1$  and not bound to collisions.
\end{nr}
\begin{proof}[Proof of~\ref{nr:temp1}]
It suffices to prove it for non-redundant symmetry groups. 
The only symmetry group of cyclic action type 
is defined by $\sigma(r) = (1,2)$ while $\rho(r)$ 
is a reflection, which is not coercive.
The same happens for brake action type:
$\rho(r)$ and $\sigma(r)$ are trivial,
and $\sigma(h) = (1,2)$. Now, $\rho(h)$ 
can be the antipodal map or the identity
(the orientation-preserving orthogonal maps  of the plane
of order less than $2$): in both cases the resulting
symmetry group is not coercive.

Now we list all the symmetry groups of dihedral type.
Up to a change in the generators, there are only
two possibilities for $\sigma$:
either $\sigma(r) = ()$ and $\sigma(h) = (1,2)$,
or $\sigma(r) = (1,2)$ and $\sigma(h) = ()$
(in fact, it cannot be $\sigma(r) = \sigma(h) = ()$
because of the transitive decomposition $2+1$, 
and the case $\sigma(r) = (1,2) = \sigma(h)$ 
can be transformed into $\sigma(r) = (1,2)$ ,
$\sigma(h')= ()$ by choosing $h'=rh$).
Let $1$, $a$ and $l$ denote respectively
the identity, the antipodal map and a reflection along a line 
in $E$. Let $l'$ denote the reflection along a line orthogonal
to $l$.
The pair $(\rho(r), \rho(h))$, by~\ref{lemma:dec}
belongs to the following list
(since the action is not of type R):
$(1,1)$, $(1,a)$, $(l,1)$, $(l,l)$, $(l,l')$, $(l,a)$,
$(a,1)$, $(a,a)$.
Hence there are $16$ possible symmetry groups.

When the action on the index set is defined by
$\sigma(r) = ()$ and $\sigma(h) = (1,2)$,
then there are the following cases:
if $(\rho(r), \rho(h))$  is in   the set 
$\{(1,1),(l,1),(l,l)\}$, then  the generated group  is bound to  collisions;
if $(\rho(r), \rho(h))$  is 
in  the set $\{ (1,a),(l,a),(l,l') \}$,
then it is  not coercive;
if $(\rho(r), \rho(h))$ is in the set 
$\{ (a,1),(a,a) \}$, then it is  redundant.
On the other hand, when $\sigma(r) = (1,2)$ and $\sigma(h) = (1,2)$,
if the pair
$(\rho(r), \rho(h))$ is in the set   
$\{ (1,1), (1,a), (l,l), (l,l'), (l,a), (a,a) \}$ then the 
action is bound to  collisions;
if $(\rho(r),\rho(h))$ is in the set  
$\{ (l,1), (a,1) \}$ then the generated group is not coercive.
This completes the proof of~\ref{nr:temp1}.
\end{proof}

At last, we now consider the symmetry groups acting transitively
on $\{1,2,3\}$.
If the action is of cyclic type, then
it has to be $\sigma(r) = (1,2,3)$.
If $\rho(r) = 1$, then 
it is the choreography symmetry group (see~\ref{sym:choreography}).
If $\rho(r) = a$, then
it is redundant (just the  choreography counted twice).
If $\rho(r) = l$, then it is the symmetry group
$C_6$ (cyclic eight -- see~\ref{sym:eights}).
Furthermore, the group cannot be of brake type, since the only groups
acting transitively and effectively on a set of three elements are 
the cyclic group of order $3$ and the  dihedral
group of order $6$.
Thus, if the action is of dihedral type, up to a permutation
in the index set one has
$\sigma(r) = (1,2,3)$ and $\sigma(h)=(1,2)$.
The only symmetry group which is of type R and 
not redundant is therefore
the Lagrange symmetry group (see~\ref{sym:lagrange}).
We want to show that the groups which are not redundant,
not bound to collisions
and not of type R are $C_6$, $D_6$ and $D_{12}$
as listed in~\ref{sym:eights}.
By~\ref{lemma:dec}, $\rho(r)$ has order at most $2$.
First consider the case $\rho(r) = a$. Since $\sigma(r)$
has order $3$, it follows that $\rho(r^3) = a$ and 
$\sigma(r^3) = ()$,
i.e.\ that the action is redundant. Thus $\rho(r)$ 
belongs to the set $\{1,l\}$,
where as above $1$ is the identity and $l$ a reflection.
If $\rho(r) = 1$, then
the element $\rho(h)$ can be $1$ or $a$.
In the first case the resulting action is bound to collisions
and in the second case it is $D_6$ (see~\ref{sym:eights}).
On the other hand, if $\rho(r) = l$, 
then $\rho(h)$ belongs to the set $\{ 1,-1,l,l'\}$,
where $l'$ is the reflection along a line orthogonal to the 
line of the reflection $l$ (it has to be an element
in $O(2)$ commuting with $l$).
Since if $\rho(h) = 1$ or $l$ then 
the action is bound to collisions,
the choices left are $-1$ and $l'$.
In both cases the resulting symmetry group is $D_{12}$ 
(see~\ref{sym:eights}). 
\end{proof}
\begin{nr}
\label{nr:nontrivialcore}
For the $3$-body problem in the plane, paths 
constrained under symmetry groups with non-trivial core,
not bound collisions and not fully uncoercive  
are always homographic.
\end{nr}
\begin{proof}
We recall that the \emph{core} of a symmetry group $G$
is $K = \ker\tau \subset G$.
The subgroup $K$ acts effectively on the index set,
hence it is isomorphic to a subgroup of $\Sigma_3$:
either its order is $2$ or $3$, or 
it is the dihedral group of order $6$ (the full
permutation group on three elements).
If its order is $2$, then it is easy to see 
that to avoid being bound to collisions or
not coercive it has to act on the plane by
the antipodal map $a$.
If its order is $3$, then necessarily its 
planar representation is given by a rotation of order $3$.
If it is a dihedral group of order $6$, then
its representation in the plane is the standard
representation of a dihedral group.
In the first case ($|K|=2$) the configuration space $\X^K$
has dimension $2$ and its elements are 
configurations $(x_1,x_2,x_3)$ where 
$x_1 = -x_2$ and $x_3=0$.
In the second case ($|K|=3$)  the configuration space $\X^K$
has dimension $2$ and its elements are
all equilateral triangles.
In the third case ($|K|=6$), the configuration space $\X^K$
has dimension $1$ and its elements are 
scalar multiples of a given equilateral triangle.
Now, as above, it is possible to add elements to $G$
which act on the time circle, but it is not difficult
to see that in the first two cases minimizers
will be, if not bound to collisions, rotating central
configurations, while in the third case the 
only way to gain coercivity is to consider 
an action bound to collisions (since it is not of type R).
\end{proof}

\begin{nr}
\label{nr:ntc}
For the $3$-body problem in the plane, if the symmetry group $G$ 
has non-trivial core, 
is not bound collisions and not fully uncoercive,  
then 
local minimizers are collisionless  
homographic solutions.
\end{nr}
\begin{proof}
As in the proof of~\ref{nr:nontrivialcore}, one is reduced to  consider only
the three cases $|K|=2,3,6$. For $|K|=6$ 
either $G$ is bound to collisions or it is fully uncoercive; 
for $|K|=2,3$
the equivariant 
$3$-body problem results to be equivalent to an equivariant
one-center planar $\alpha$-homogeneous 
Kepler problem, and the thesis follows since local minima
of the equivariant Kepler problem are collisionless (see 
remark~\ref{rem:kepler} below).
\end{proof}
\begin{remark}\label{rem:kepler}
Consider the one-center Kepler problem in the plane and 
a symmetry group $G$ of the action functional. In this 
case $\X = E \cong \R^2$, and obviously the action on the
index set $\{1\}$ is trivial.
If $G$ has non-trivial core
then $\X^{\ker\tau}$ has dimension $0$ or $1$, and hence 
either $G$ is bound to collisions or fully uncoercive. 
Thus, $\ker\tau=1$ and 
so  in order to show that local minimizers are collisionless, 
one needs  only  to exclude boundary collisions (see~\ref{T1}
below).
Let  $\rho(h)$ be the image
under $\rho$ of  the generator $h$ of the $\T$-isotropy in question.
If $\rho(h)$ is a rotation of angle $\pi$, 
then $G$ is bound to collisions. Otherwise,
$\rho(h)$ is a reflection along a line, and therefore 
the claim follows by  theorem~\ref{T2} in section~\ref{sec:sec7}.
\end{remark}

As a consequence of theorem~\ref{MT1}, one can easily proof the
following proposition.
\begin{lemma}
\label{lemma:lag}
Let $G$ be a symmetry group such that $\action_\omega^G$ is coercive.
If  $\Lambda^G$ contains a Lagrange relative equilibrium motion $x$,
then the minimum of the action functional $\action_\omega^G$ 
is attained on $x$ (and hence the minimizer is homographic).
\end{lemma}

\begin{remark}
It follows from~\ref{le:LS},~\ref{le:21C} and~\ref{lemma:nottypeR}
that a sort of converse of~\ref{lemma:lag} holds, at least 
in the case of almost equal masses:
That is, for every $G$ 
if the masses are approximately equal (subject
to the symmetry constraints) and  $\Lambda^G$ 
does not contain Lagrange equilibrium motions, then there always
exists a suitable choice of angular speed $\omega$ such 
that the minimum of $\action^G_\omega$ is not homographic --
that is, it is not a rotating collinear Euler configuration.
It is of some interest to understand whether this is true for 
all values of masses.
\end{remark}

\section{Parabolic collisions with isosceles symmetry}
\label{sec:parab}

In this section we will study a special class of colliding trajectories,
constrained under a simple symmetry.  This symmetry, which we term
\emph{isosceles}, is the only one not fulfilling the rotating circle property.
Our aim is to prove that also for this symmetry there exists a local variation
which implies that minimizers are collisionless.  
In the following we can
assume  $\alpha \in (0,2)$,
$\kk \supset \{1,2\}$ and  $m_1=m_2$.
\begin{definition}\label{parab:coll:traj} 
A \emph{parabolic collision trajectory} 
for the cluster $\kk\subseteq \nn$
is the path  
\begin{equation*} 
{q}_i(t) = |t|^{2/(2+\alpha)} \xi_i, \quad i \in \kk,\,\, t \in \R 
\end{equation*} 
where $\xi=(\xi_i)_{i \in \mathbf{k}}$ 
is a central configuration with $k$ bodies.  
\end{definition} 

\begin{definition}\label{escaping:path} 
An \emph{escaping path} for the cluster $\kk\subseteq \nn$
is a path of type $y=q+\varphi$, 
where $q$ is a parabolic collision trajectory for the cluster $\kk$
and $\varphi \in H_0^1(\R)$.
\end{definition} 

We now define the action of $g_0$ on the time line $\R$,
on the space $E$ and on the index set as follows:
$g_0(t) = -t $  for $t\in \T$, 
$g_0$ 
acts as a reflection along a line in $E$ and  
$\sigma(g_0) = (1,2)$.

\begin{definition}\label{escaping:equiv:path} 
Given $g_0$, we say that an \emph{escaping path} is
\emph{$g_0$-equivariant} if 
\begin{equation*} 
y(g_0t) = g_0 y(t).
\end{equation*} 
\end{definition} 

Notice that Definition~\ref{escaping:equiv:path} 
covers both the cases of
binary and triple clusters ($k=2,3$).

\begin{remark}
Let $l$ the line fixed by $\rho(g_0)$; then $y_1(0)$ and $y_2(0)$ are
symmetric with respect to $l$ while $y_3(0)$ (if $n>2$) belongs to $l$.  A
{$g_0$-equivariant} {escaping path} is determined by its restriction on 
the half line $[0,+\infty)$ provided that $y(0) \in \mathcal{X}^{g_0}$.
\end{remark}

We term $\mathcal{L}_{\mathbf{k}}$ the partial Lagrangian function, 
when $x=(x_i)_{i \in \mathbf{k}}$, $\mathcal{L}_{\mathbf{k}}(x)$ is 
the Lagrangian restricted on the bodies of the cluster ${\mathbf{k}}$.

\begin{definition}
\label{min_par_coll_traj}
We say that a parabolic collision trajectory, 
$q=(q_i)_{i\in \mathbf{k}}$, defined in~\ref{parab:coll:traj}, 
is a \emph{$g_0$-equivariant minimizing parabolic collision trajectory}
if for every $g_0$-equivariant escaping path 
$y=q+\varphi$ the integral of the variation is positive:
\[
\int_{-\infty}^{+\infty}[ 
\mathcal{L}_{\mathbf{k}}(q+\varphi)- \mathcal{L}_{\mathbf{k}}(q)] dt
\geq 0.
\]
\end{definition} 

Let $\delta \in (\R^2)^k$, $k=2,3$, be a vector of norm 
$|\delta| = (\sum_{i=1}^k \delta_i^2 )^{1/2}$ 
sufficiently small and $T>0$ a real
number. 

\begin{definition}
\label{stand:var}
The \emph{standard variation}
associated to $\delta$ and $T$ is defined as 
\begin{equation*}
v^\delta(t)= \left\{
\begin{array}{ll}
\delta & \mbox{if} \,\,\, 0 \leq |t| \leq T-|\delta| \\
(T-t)\dfrac{\delta}{|\delta|} & \mbox{if} \,\,\, T-|\delta| \leq |t| \leq T \\
0 & \mbox{if} \,\,\, |t| \geq T. \\
\end{array}
\right.
\end{equation*}
\end{definition}

\begin{remark}
Let $q$ be a $g_0$-equivariant minimizing parabolic collision trajectory 
and $v^\delta(t)$ a standard variation. 
Then the path $y(t)=q(t)+v^\delta(t)$ (with $t \in \R$)
is $g_0$-equivariant if and only if  $v^\delta$ is $g_0$-equivariant,
and hence  $\delta=g_0\delta$.
Thus, in particular, any $\delta$ fixed by $g_0$ yields a 
$g_0$-equivariant standard variation $v^\delta(t)$.
\end{remark} 

In the next theorem we will prove 
that $g_0$-equivariant parabolic
collision trajectories 
cannot be local minimizers
(see~definition \ref{min_par_coll_traj}).

\begin{theorem}
\label{main_thm}
Let $q$ be a $g_0$-equivariant parabolic collision trajectory.
Then there exists a $g_0$-equivariant standard variation $v^\delta$ such that
the path $q+v^\delta$ does not have a collision at $t=0$ and
\[
\Delta \mathcal{A} = 
\int_{-\infty}^{+\infty} [\mathcal{L}_{\mathbf{k}}(q+v^\delta) - 
\mathcal{L}_{\mathbf{k}}(q)] dt < 0.
\]
\end{theorem}

The proof of this result requires several intermediate steps.
To begin with, consider the function
\begin{equation}
S(\xi,\delta)=\int_0^{+\infty} 
\left( \frac{1} {\left| \xi t^{2/(2+\alpha)}-\delta \right|^{\alpha}} -
\frac{1}{\left| \xi t^{2/(2+\alpha)} \right|^{\alpha}}  \right) dt
\end{equation}
where $\xi,\delta \in \R^2$.

The next result allows to estimate the action differential 
involved in a standard variation.

\begin{prop}
\label{action_var}
Let $q=\{q\}_i= \{t^{2/(2+\alpha)} \xi_i \}$, 
$i=1,\ldots,k$ be a parabolic collision 
trajectory and $v^\delta$ a $g_0$-equivariant standard variation. 
Then for $\delta \rightarrow 0$
\begin{equation*}
\Delta {\cal A} = 2|\delta|^{1-\alpha/2} 
\sum_{\substack{i<j \\ i,j \in \kk}} 
m_i m_j S(\xi_i -\xi_j, \frac{\delta_i - \delta_j}{|\delta|}) + O(|\delta|).
\end{equation*}
\end{prop}

To prove Theorem~\ref{main_thm} 
we have to provide a suitable 
$g_0$-equivariant
standard variation such that the right end side of~\ref{action_var} is
negative. This will depend 
on the particular central configuration $\xi$ drawn
by the parabolic collision trajectory. 
The proof of~\ref{action_var}
can be found 
in  \cite{FT2003}, Lemma (9.2). 
Let $\vartheta \in [0,2\pi]$ 
such that 
$\cos \vartheta=\left\langle 
{\xi}/{\left|\xi\right|},\delta\right\rangle$ then 
\begin{equation}
\label{SpropPhi}
S(\xi,\delta)=
\left|\xi\right|^{-1-\alpha/2}\int_0^{+\infty} 
\frac{1}{\left( t^{\frac{4}{\alpha+2}} 
-2\cos \vartheta t^{\frac{2}{\alpha+2}} 
+1 \right)^{\alpha /2}} - 
\frac{1}{t^{\frac{2\alpha}{\alpha+2}}} dt.
\end{equation}

\begin{prop}
The function
\begin{equation}
\label{Phi1}
\Phi_\alpha(\vartheta)=\int_0^{+\infty} 
\frac{1}{\left( t^{\frac{4}{\alpha+2}} 
-2\cos \vartheta t^{\frac{2}{\alpha+2}} 
+1 \right)^{\alpha /2}} - 
\frac{1}{t^{\frac{2\alpha}{\alpha+2}}} dt, 
\quad \alpha \in (0,2)
\end{equation}
is defined and continuous 
on the interval $(0,2\pi)$ and satisfies the following properties
\begin{enumerate}
\item%
\label{p:i}
$\Phi_\alpha(\vartheta)=\Phi_\alpha(2\pi-\vartheta)$, 
for every $\vartheta \in (0,2\pi)$ and $\alpha \in (0,2)$, 
i.e.\ its plot is symmetric with respect to $\vartheta=\pi$;
\item%
\label{p:ii}
is decreasing on $(0,\pi]$ and 
increasing on $[\pi,2\pi)$, i.e.\ it 
achieves its minimal value on $\vartheta = \pi$;
\item%
\label{p:iii}
when $\alpha < 1$ then $\Phi_\alpha(0)=\Phi_\alpha(2\pi)$
is finite, while when $\alpha \geq 1$ 
\[
\lim_{\vartheta \rightarrow 0^+}\Phi_\alpha(\vartheta)
=\lim_{\vartheta \rightarrow 2\pi^-}\Phi_\alpha(\vartheta)=+\infty.
\]
\end{enumerate}
\end{prop}
\begin{proof}
Properties~\ref{p:i} and~\ref{p:ii}
are obvious.
To prove~\ref{p:iii}, we study the integral of the function 
\[
\phi_\alpha(t,\vartheta)= 
\frac{1}{\left( t^{\frac{4}{\alpha+2}} 
-2\cos \vartheta t^{\frac{2}{\alpha+2}} +1 \right)^{\alpha /2}} 
- \frac{1}{ t^{\frac{2\alpha}{\alpha+2}}}.
\]
on the interval $(0,+\infty)$. For 
$\vartheta \in (0,2\pi)$,
the function $\phi_\alpha$ is a continuous with respect to $t$ in
$(0,+\infty)$ and we have to take into
account the convergence of its integral at 0 and at $+\infty$.

When $t \rightarrow 0^+$, the integral of 
$\phi_\alpha(\cdot,\vartheta)$ 
in a right neighborhood of 0 is finite, being 
${2\alpha}/({2+\alpha})<1$, 
$\forall \alpha \in (0,2)$. 
About the integrability at infinity,  we have 
\[
\phi_\alpha(t,\vartheta)= 
\frac{1}{ t^{\frac{2\alpha}{\alpha+2}}}\left[ 
\frac{1}{\left( 1 -2\cos \vartheta t^{\frac{-2}{\alpha+2}} + 
t^{\frac{-4}{\alpha+2}} \right)^{\alpha /2}} - 1 
\right] 
\approx 
\frac{\alpha}{2}\left( 
\frac{1}{t^2}-2\frac{\cos \vartheta}{t^{2\frac{\alpha+1}{\alpha+2}}} 
\right)
\]
as $t \rightarrow +\infty$; since
for every $\alpha \in (0,2)$, 
$2{(\alpha+1)}/{(\alpha+2)} >1$ and 
$\Phi_\alpha$ is integrable on $[2,+\infty)$.

When $\vartheta =0$, the function 
$\phi_\alpha(\cdot,0)=\phi_\alpha(\cdot,2\pi)$ has a
singularity of order $\alpha$ in $t=1$. Therefore the integral function
$\Phi_\alpha$ 
tends to $+\infty$ as $\vartheta$ tends to $0$ 
if and only if $\alpha \geq 1$.
\end{proof}

\begin{corollary}
There exists 
$\bar{\vartheta}=\bar{\vartheta}(\alpha) \in (0,\pi/2)$ 
such that
$\Phi_\alpha(\bar{\vartheta})=0$ and 
$\Phi_\alpha(\vartheta)<0$, for every
$\bar{\vartheta}<\vartheta<2\pi-\bar{\vartheta}$.
\end{corollary}

\subsection*{Binary collisions}
The simplest case is when only two bodies are involved 
in the collision; thus $\xi_1 = -\xi_2$ and we have to take
into account the angle between $\xi_1$ and the line $l$ fixed
by $g_0$.
\begin{figure}[ht]
\psfrag{xi1}{$\xi_1$}
\psfrag{xi2}{$\xi_2$}
\psfrag{de}{$\delta$}
\psfrag{md}{$-\delta$}
\psfrag{th}{$\vartheta$}
\begin{center}
\includegraphics[width=0.5\textwidth]{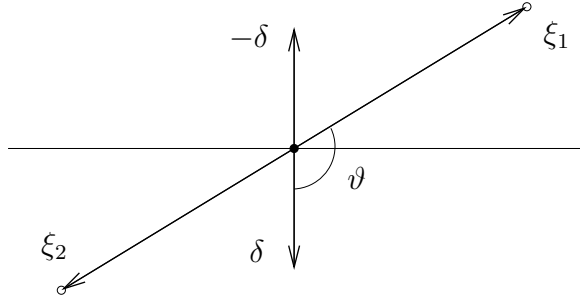}
\caption{Binary collisions}
\label{fig:binary:coll}
\end{center}
\end{figure}
Consider a vector $\delta$ orthogonal 
to $l$ with norm $|\delta|={1}/{2}$, in such a way that the
angle $\vartheta = \arccos \left\langle \frac{\xi_1 -\xi_2}{\left| \xi_1
-\xi_2 \right|}, \delta\right\rangle$ lies in the interval
$[{\pi}/{2},\pi]$ (see Figure~\ref{fig:binary:coll}). 
To obtain a negative  variation of the
interaction energy, we add to the position of the first body a vector
$-\delta$ and to the second a vector $\delta$.

\subsection*{Triple collisions}
A triple collision can take place in two different ways: from collinear
configurations or from the Lagrange configuration. We study separately these
two cases.  Our aim is to show that is always possible to reduce a triple
collision to a binary one.

\subsubsection*{Triple collisions from collinear central configurations}
\begin{figure}[ht]
\begin{center}
\subfigure[]{%
\psfrag{th}{$\vartheta$}
\psfrag{de}{$\delta$}
\psfrag{xi1}{$\xi_1$}
\psfrag{xi2}{$\xi_2$}
\psfrag{xi3}{$\xi_3$}
\psfrag{l}{$l$}
\includegraphics[width=0.4\textwidth]{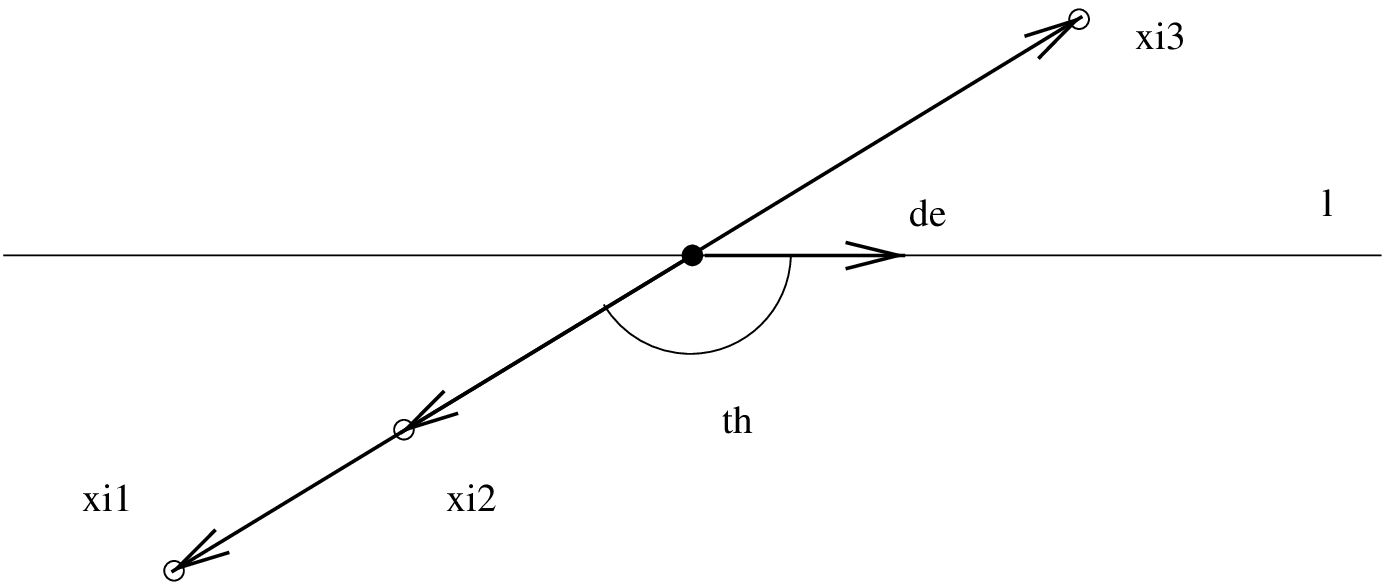}
\label{subf:a}
}
\subfigure[]{%
\psfrag{th}{$\vartheta$}
\psfrag{de}{$\delta$}
\psfrag{xi1}{$\xi_1$}
\psfrag{xi2}{$\xi_2$}
\psfrag{xi3}{$\xi_3$}
\psfrag{md}{$-\delta$}
\psfrag{l}{$l$}
\includegraphics[width=0.4\textwidth]{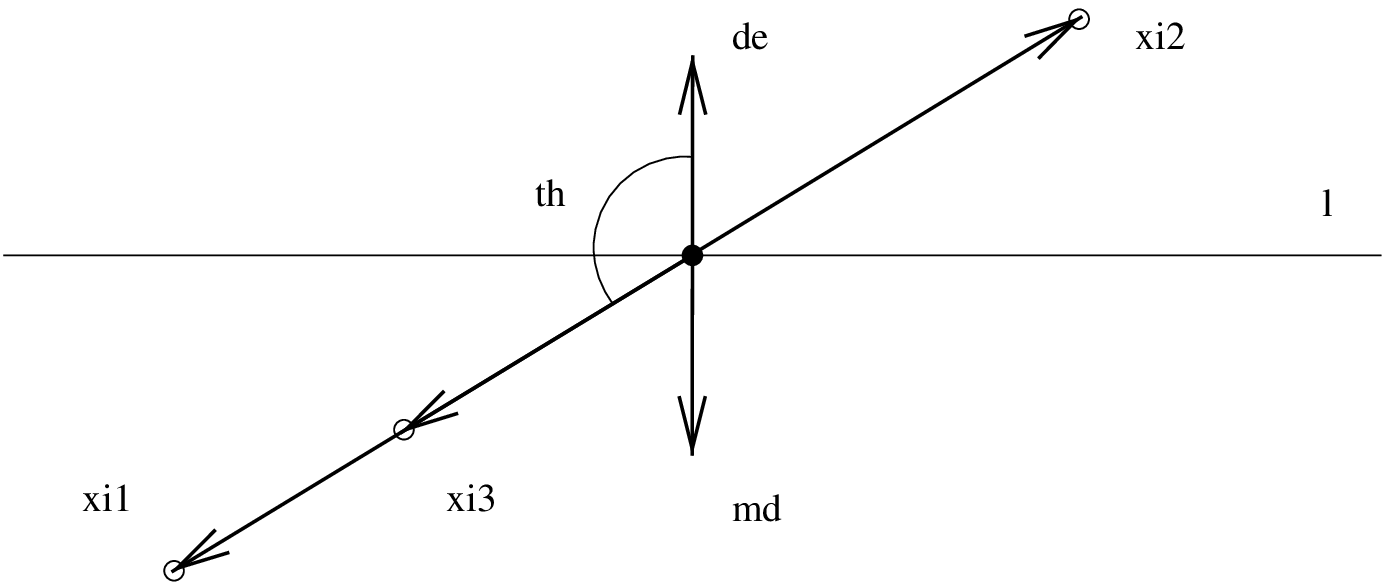}
\label{subf:b}
}
\caption{Triple collision from collinear central configuration}
\label{fig:triple}
\end{center}
\end{figure}

Two different situations occur 
whether the second or the third body lies
between the other two (respectively we will have
$\xi_1=-\lambda\xi_3$ or $\xi_1=-\lambda\xi_2$, $\lambda>0$).  
In the first
case we refer to Figure~\ref{subf:a} 
and we move the trajectory of the third body of
a vector $\delta$ parallel 
to $l$ whether $\xi_2=\mu \xi_1$, or $\xi_2=\mu \xi_3$,
$\mu \in [0,1)$. 

When the third body lies in the middle, let $\delta$ 
be orthogonal to  $l$,
$|\delta|={1}/{2}$, such that $\vartheta = \arccos \left\langle
\frac{\xi_1 -\xi_2}{\left| \xi_1 -\xi_2 \right|}, \delta\right\rangle$ is in
the interval $[{\pi}/{2},\pi]$, see 
Figure~\ref{subf:b}. If we shift the second
body of $\delta$ and the first  of $-\delta$, the interaction between these
two bodies decreases. Since $\vartheta > \pi/2$ and $\mu<1$, also the sum of
the
variations of the interactions of the third body with the other two is
negative.
\subsubsection*{Triple collisions from triangular central configurations}
\begin{figure}[hb]
\begin{center}
\psfrag{ga}{$\gamma$}
\psfrag{de}{$\delta$}
\psfrag{d1}{$\xi_1-\xi_3$}
\psfrag{d2}{$\xi_2-\xi_3$}
\psfrag{h}{$h$}
\psfrag{l}{$l$}
\includegraphics[width=0.5\textwidth]{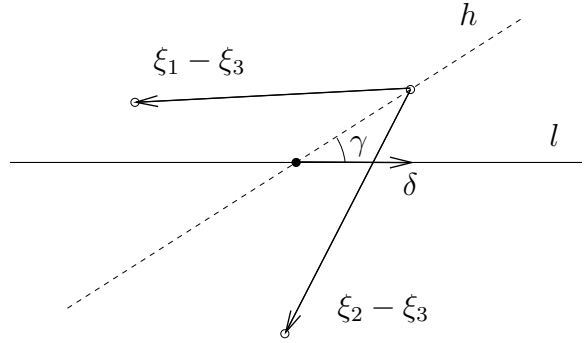}
\caption{Triple collisions from Lagrange configuration.}
\label{fig:triple:lagrange}
\end{center}
\end{figure}

Consider now the $3$ particles at the vertices of a regular triangle moving on a
parabolic collision trajectory that makes them collide at $t=0$ in their
center of mass. Since $m_1=m_2$, the center of mass of the cluster lies on the
perpendicular line, $h$, from the third body to the line joining the other
two.  Our aim is to show that there always exists a vector $\delta \in S^1$,
such that, when  replacing $\xi_3$ with $\xi_3 + \delta$,
the interaction 
potential decreases (the
sum of the two variations is negative even thought they are
not necessarily both
negative) and therefore $\Delta\action$ 
of~\ref{action_var} is negative. In the following
we refer to figure~\ref{fig:triple:lagrange}
and we call $\gamma$ the
angle in $[0,\pi/2]$ with edges the line $h$ and the subspace $l$. We will
prove the following result.

\begin{propo}
\label{thm3coll}
For every $\alpha \in (0,2)$ the following inequality holds
\begin{equation}
\label{dis:thm}
\Phi_\alpha(\frac{2\pi}{3}+\gamma)+
\Phi_\alpha(\frac{2\pi}{3}-\gamma)<0, \quad \forall \gamma \in \left[0,\frac{\pi}{2}\right].
\end{equation}
\end{propo}

Before proving Proposition~\ref{thm3coll} we need two preliminary Lemmata. 

\begin{lemma}
\label{3coll:le1}
For  every $\alpha \in (0,2)$ and $\vartheta \in (0,2\pi)$
\begin{multline}
\label{ug:le1} 
\frac{2}{\alpha(\alpha+2)}
\Phi_\alpha(\vartheta)  = \frac{1}{\alpha-2}          
\beta(\frac{\alpha+2}{4},\frac{\alpha+2}{4}) +{}\\{}+ 
\frac{1}{\alpha} \sum_{k=1}^{+\infty} 
\binom{-\alpha/2}{k}
(-1)^k 2^k (\cos \vartheta)^k 
\frac{\alpha +2k}{\alpha+2k-2}
\beta(\frac{\alpha+2}{4}+\frac{k}{2},\frac{\alpha+2}{4}
+\frac{k}{2}),
\end{multline}
where the function $\beta : \R \times \R \rightarrow \R^+$ is defined as 
\begin{equation}
\label{beta}
\beta(z,w) = \int_0^1 t^{z-1}(1-t)^{w-1}dt.
\end{equation}
\end{lemma}

\begin{remark}
\label{prop:beta}
It is well known that 
the function $\beta$ satisfies the following properties
\begin{enumerate}
\item\label{en:a}%
$\beta(z,w)=\beta(w,z)$, for every $z,w \in \R$;
\item\label{en:b}
$\beta(z+1,w)=\frac{z}{z+w}\beta(z,w)$, for every $z,w \in \R$;
\item\label{en:c}%
$\beta(x+n,x+n)=\frac{1}{4^n}\frac{
\binom{-x}{n}
}
{%
\binom{-x-1/2}{n}
}%
\beta(x,x)$, for every $x \in \R$, $n \in \N$.
\end{enumerate}
\end{remark}

\begin{proof}[Proof of Lemma~\ref{3coll:le1}]
Consider the function $\Phi_\alpha$ as in~\ref{Phi1} with
$t=s^{\frac{2+\alpha}{2-\alpha}}$: then
\begin{eqnarray*}
\Phi_\alpha(\vartheta) & = & \frac{2+\alpha}{2-\alpha} \int_0^{+\infty} 
[ \frac{s^{\frac{2\alpha}{2-\alpha}}}
{\left( s^{\frac{4}{2-\alpha}} -2 s^{\frac{2}{2-\alpha}} \cos \vartheta 
+1 \right)^{\alpha / 2}}- 1 ] ds \\
& = &  \frac{2+\alpha}{2-\alpha} \int_0^{+\infty} 
[ \frac{s^{\frac{2\alpha}{2-\alpha}}}
{( s^{\frac{4}{2-\alpha}}+1)^{\alpha/2}  
(1- \frac{2 s^{\frac{2}{2-\alpha}} 
\cos \vartheta}{s^{\frac{4}{2-\alpha}}+1} 
)^{\alpha/2}}- 1 ] ds.
\end{eqnarray*}
Since 
$x\leq \frac{1}{2}(x^2+1)$ and, for $|z|<1$ 
there holds
$(1+z)^{\alpha} = 
\sum_{k=0}^{+\infty}
\binom{\alpha}{k} z^k$,
we can write
\begin{equation*}
\Phi_\alpha(\vartheta)   =   
\frac{2+\alpha}{2-\alpha} 
\int_0^{+\infty} 
\left\{ 
 \frac{s^{\frac{2\alpha}{2-\alpha}}}%
{(s^{\frac{4}{2-\alpha}}+1 )^{\alpha/2}}  +
\sum_{k=0}^{+\infty} 
\binom{-\alpha/2}{k}
\frac{%
(-1)^k 2^k 
s^{\frac{2k}{2-\alpha}}
(\cos \vartheta)^k
} {( s^{\frac{4}{2-\alpha}}+1)^k} 
-1\right\} ds  
\end{equation*}
or, equivalently,
\begin{equation}
\label{Phi2} 
\Phi_\alpha(\vartheta)   =   
\frac{2+\alpha}{2-\alpha} \left\{ I_0 + 
\sum_{k=1}^{+\infty} 
\binom{\alpha/2}{k}
(-1)^k 2^k  (\cos \vartheta)^k  I_k\right\}, 
\end{equation}
where
\begin{eqnarray}
\label{I_0}
I_0 &=& \int_0^{+\infty} 
\left[ \frac{s^{\frac{2\alpha}{2-\alpha}}}%
{( s^{\frac{4}{2-\alpha}}+1 )^{\alpha/2}}-1 \right]ds, 
\label{I_k}\\
I_k &=& \int_0^{+\infty} \frac{s^{\frac{2\alpha}{2-\alpha}}}%
{( s^{\frac{4}{2-\alpha}}+1)^{\alpha/2}}%
\frac{s^{\frac{2k}{2-\alpha}}}{( s^{\frac{4}{2-\alpha}}+1)^k}%
ds, \quad k \in \N^*.
\end{eqnarray}
We now compute 
the expressions of $I_0$ and $I_k$ separately to 
obtain~\ref{ug:le1}. Integrating by parts,~\ref{I_0} becomes
\begin{equation*}
I_0 = -\frac{2\alpha}{2-\alpha}
\int_0^{+\infty}\frac{s^{\frac{2\alpha}{2-\alpha}}}
{(s^{\frac{4}{2-\alpha}}+1)^{\frac{\alpha}{2}+1}}ds;
\end{equation*}
we then substitute
\begin{equation}
\label{subst}
s=( \frac{1-t}{t})^{\frac{2-\alpha}{4}}, \quad 
\frac{ds}{dt}=\frac{2-\alpha}{4}
( \frac{t}{1-t} )^{\frac{2+\alpha}{4}}
( -\frac{1}{t^2} )
\end{equation}
to obtain
\begin{equation}
I_0 = -\frac{\alpha}{2}\int_0^{1} 
t^{\frac{\alpha-2}{4}} (1-t)^{\frac{\alpha-2}{4}} dt =
-\frac{\alpha}{2} \beta 
( \frac{\alpha+2}{4},\frac{\alpha+2}{4} ). 
\label{I_0n}
\end{equation}
From~\ref{I_k}, by~\ref{subst}, one deduces that 
\begin{eqnarray}
\label{I_kn}
I_k & = & \frac{2-\alpha}{4} 
\int_0^1 t^{\frac{\alpha}{4}+\frac{k}{2}-\frac{3}{2}} 
(1-t)^{\frac{\alpha}{4}+\frac{k}{2}-\frac{1}{2}} dt \notag \\
& = & \frac{2-\alpha}{4} \beta 
( \frac{\alpha}{4}- \frac{1}{2} + \frac{k}{2},%
\frac{\alpha}{4}+ \frac{1}{2} + \frac{k}{2} ). 
\end{eqnarray}
By replacing~\ref{I_0n} and~\ref{I_kn} in~\ref{Phi2} 
it follows that 
\begin{multline*}
\Phi_\alpha(\vartheta) = \frac{\alpha(\alpha+2)}{2}
\left\{ \frac{1}{\alpha-2}
\beta\left(\frac{\alpha+2}{4},\frac{\alpha+2}{4}\right) 
\right.{}+\\{}+ \left.  
 \frac{1}{\alpha} \sum_{k=1}^{+\infty}
\binom{-\alpha/2}{k}
(-1)^k 2^{k-1} (\cos \vartheta)^k \beta\left( 
\frac{\alpha}{4}- \frac{1}{2} + \frac{k}{2},%
\frac{\alpha}{4}+ \frac{1}{2} + \frac{k}{2} \right) \right\}. 
\end{multline*}
Since $\beta(z-1,z)=\frac{2z-1}{z-1}\beta(z,z)$ 
(see~\ref{en:b} in Remark~\ref{prop:beta}), the claim~\ref{ug:le1} 
is proved.
\end{proof}

\begin{lemma}
\label{3coll:le2}
For every $x \in (1/2,1)$, the following inequality holds
\begin{equation}
\label{dis:le2}
\frac{1}{x-1} + \sum_{k=1}^{+\infty} 
{\binom{-x}{k}}^2
\Big(\frac{3}{4}\Big)^k \frac{1}{x+k-1}\frac{4^k(k!)^2}{(2k)!}<0.
\end{equation}
\end{lemma}
\begin{proof}
Define the function $f_k$ as  $f_k(x)=
{\binom{-x}{k}}^2
\frac{1}{x+k-1}\frac{4^k(k!)^2}{(2k)!}$ for every
$x \in (1/2,1)$ and  $k \in \N$. 
For a fixed $x \in (1/2,1)$, $f_k(x)$ 
is monotonically decreasing in $k$, while, for a fixed $k$,
$f_k(x)$ is  increasing  in $x$ on $(1/2,1)$.
We first prove that the following inequality holds
\begin{equation} 
\label{serie:troncata} 
\sum_{k=5}^{+\infty}f_k(x)\Big(\frac{3}{4}\Big)^k < \dfrac{27}{35}, \qquad \text{for all } x \in (1/2,1).
\end{equation}
Indeed, 
by monotonicity of $f_k(x)$, 
inequality~\ref{serie:troncata} is implied by the fact that 
\[
f_5(1) \sum_{k=5}^{+\infty}
\Big(\frac{3}{4}\Big)^k = f_5(1) \Big(\frac{3}{4}\Big)^5 4
=\frac{4}{5} \frac{3^5 (5!)^2}{(10)!} = \dfrac{27}{35} 
\] 
Hence~\ref{dis:le2}  would follow once it is proved that,
for every $x\in[1/2,1]$
\begin{equation}
\label{pol:est}
\frac{1}{x-1} + 
\sum_{k=1}^{4}f_k(x) \Big(\frac{3}{4}\Big)^k +\dfrac{27}{35}<0.
\end{equation}
Expanding the expression above we can write
\begin{multline*}
\dfrac{1}{x-1} + \frac{3}{2}x + \frac{3}{8}(x+1)x^2
+\frac{3}{80}(x+2)(x+1)^2x^2 +{}\\ {}+
\frac{9}{4480}(x+3)(x+2)^2(x+1)^2x^2 + \dfrac{27}{35} < 0
\end{multline*}
or, equivalently,~\ref{dis:le2} holds if for every $x\in (1/2,1)$
\begin{multline}\label{eq:ineq}
p(x) = 9x^8 + 72x^7 + 
366x^6 + 
684 x^5 +{} \\ {}+
1749 x^4 -
756 x^3 +
4596 x^2 -
3264 x + 1024 
> 0;
\end{multline}
one way of proving \ref{eq:ineq} is just to show
that the polynomial $p(x)$ has a Taylor expansion
centered in $1/2$
with positive coefficients; alternatively, 
since in $(1/2,1)$ 
\begin{equation*}
p(x) \geq 1024 - 3264x +  
4596 x^2 -
756 x^3  
\geq  1024 \left(
1 - \frac{10}{3}x + 4x^2 - x^3
\right),
\end{equation*}
the claim follows from the fact that 
the minimum of the cubic polynomial $1-\frac{10}{3}x + 4x^2 - x^3$ 
in $(1/2,1)$ is attained at $x_0=\frac{4-\sqrt{6}}{3}$,
with value
$\frac{35}{27} - \frac{4}{9}\sqrt{6}>0$.
\end{proof}

\begin{lemma}
\label{le:pi/6}
For every $\alpha \in (0,2)$ the following inequality holds
\begin{equation}
\label{dis:pi/6}
\Phi_\alpha\Big(\frac{\pi}{6}\Big)+\Phi_\alpha\Big(\frac{7\pi}{6}\Big) < 0.
\end{equation}
\end{lemma}

\begin{proof}
We replace $\gamma = \frac{\pi}{2}$ in~\ref{dis:thm} and we use the result
proved in Lemma~\ref{3coll:le1} to obtain
\begin{multline}
\Phi_\alpha\Big(\frac{\pi}{6}\Big)+\Phi_\alpha\Big(\frac{7\pi}{6}\Big) =
 \frac{1}{\alpha-2}\beta\left( \frac{\alpha+2}{4},\frac{\alpha+2}{4}\right) 
+{}\\{}+ 
 \frac{1}{\alpha} \sum_{k=1}^{+\infty} 
\binom{-\alpha/2}{2k}
\Big( \frac{3}{4}\Big)^{k} 2^{2k} 
 \frac{\alpha +4k}{\alpha +4k -2}
\beta\left(\frac{\alpha+2}{4}+k,\frac{\alpha+2}{4}+k\right). 
\end{multline}
By~\ref{en:c} of Remark~\ref{prop:beta} one obtains
\[
\beta\left(\frac{\alpha+2}{4}+k,\frac{\alpha+2}{4}+k\right)=
\frac{1}{4^k}\frac{%
\binom{-\frac{\alpha+2}{4}}{k}
}{%
\binom{-\frac{\alpha+2}{4}-1/2}{k}
} 
\beta(\frac{\alpha+2}{4},\frac{\alpha+2}{4})
\]
and then~\ref{dis:pi/6} is implied by
\begin{equation*}
\frac{1}{\alpha-2} + 
\frac{1}{\alpha}\sum_{k=1}^{+\infty} 
\binom{-\alpha/2}{2k}
\Big( \frac{3}{4}\Big)^{k} \frac{\alpha+4k}{\alpha+4k-2} 
\frac{%
\binom{-\frac{\alpha+2}{4}}{k}
}{%
\binom{-\frac{\alpha+2}{4}-1/2}{k}
} < 0.
\end{equation*}
Finally, this inequality is equivalent to 
\begin{equation*}
\frac{1}{\alpha-2} + \sum_{k=1}^{+\infty} \Big( \frac{3}{4}\Big)^{k} 
{\binom{-\frac{\alpha+2}{4}}{k}}^2
\frac{1}{\alpha+4k-2}\frac{4^k(k!)^2}{(2k)!} < 0.
\end{equation*}
To conclude the proof we  use lemma~\ref{3coll:le2} with  $x=({\alpha +2})/{4}$.
\end{proof}
\begin{proof}[Proof of Proposition~\ref{thm3coll}]
At first we remark that,
when $\gamma \in [0,\pi/6]$, then 
${2\pi}/{3} \pm \gamma \in [\pi/2,\pi]$ and therefore,
since both $\Phi_\alpha({2\pi}/{3} \pm \gamma)$ are negative, the thesis follows. 
Well then, let $\gamma \in (\pi/6,\pi/2]$ and
\begin{equation*}
\delta=\frac{2\pi}{3}-\gamma, \qquad \delta'=\frac{2\pi}{3}+\gamma
\end{equation*}
Then $\delta \in [\pi/6,\pi/2)$ and $\delta' \in (5\pi/6,7\pi/6]$. Since $\Phi_\alpha$
is decreasing on $(0,\pi)$ and is symmetric with respect to $\vartheta=\pi$ we have that
\begin{equation*}
\Phi_\alpha(\delta)\leq \Phi_\alpha\left(\frac{\pi}{6}\right), \; \forall \delta \in \left[\frac{\pi}{6},\frac{\pi}{2}\right)
\qquad \text{and} \qquad
\Phi_\alpha(\delta')\leq\Phi_\alpha\left(\frac{7\pi}{6}\right), \; \forall \delta' \in
\left(\frac{5\pi}{6},\frac{7\pi}{6}\right].
\end{equation*}
We conclude using the result we proved in Lemma~\ref{le:pi/6}, indeed
\[
\Phi_\alpha(\delta)+\Phi_\alpha(\delta') \leq \Phi_\alpha(\pi/6)+\Phi_\alpha(7\pi/6)<0.
\]

\end{proof}

\section{Proof of Theorem~\ref{MT2}}
\label{sec:sec7}

We have already described (Proposition~\ref{prop:2.27}) a criterion that
guarantees 
the existence of minimizers, we now rule out the occurrence of collisions. 
Following \cite{FT2003}, we introduce 
the key property of the action of a  finite group 
$G$ on the loop space $\Lambda$ (we give here the 
definition specialized for the $3$-body problem in the plane). 

\begin{definition}  
\label{rcp}
We say that a finite group $K$ 
acts  on $\nn = \{ 1,2,3 \}$  (resp.\ $\{1,2\}$) and $E$ with
the \emph{rotating circle property (RCP)} if 
the following conditions are true:
\begin{enumerate}
\item
for every $g \in K$  the determinant 
$\det(\rho(g)) = 1$;
\item\label{proper:2}
there exist at least two different indices
$i_1,i_2 \in \nn$  such that 
$\forall g\in K, 
(gi_1=i_1 \vee gi_2=i_2) \implies
\rho(g) = 1$ (resp.\ there exists at least 
one index $i_1 \in \{1,2\}$ 
such that  $\forall g\in K: gi_1=i_1 \implies \rho(g) = 1$).
\end{enumerate}
\end{definition}

\begin{remark}
A group $G$ has the RCP if all its $\T$-isotropy
subgroups have the RCP, 
and a trivial group has the RCP.
The underlying idea is the following:
given a trajectory with a colliding cluster $\kk \subset \nn$,
since at least two indices in $\nn$ satisfy property~\ref{proper:2}, 
after a blow-up one is left 
with a parabolic $K$-equivariant collision-ejection 
trajectory, where $K$ is a suitable subgroup of $G$.
Now, if $G$ has the RCP, then it follows that $K$ has the RCP.
The first attempt would be  to move away one of
the particles, say $x_{i_1}$, and to let it rotate in a circle.
Thanks to  the averaging estimate \cite{FT2003} 
it would be possible to prove
that the action decreases. But, this would yield 
a non-symmetric path, in general. So, it is necessary to move
at the same time all those particles $x_{gi_1}$,
for $g\in K$, which are related
by the $K$-action to the $i_1$-particle.
\end{remark}

\begin{definition} \label{defi:interiorboundarycollision} 
We say that $x \in \Lambda^G$ has an 
\emph{interior collision at $t$}
if $t \in (\I\smallsetminus \partial \I) \cap
\Delta^{-1}x$.  We say that $x$ has an 
\emph{boundary collision at $t$} 
if $t \in \partial \I\cap \Delta^{-1}x$. 
\end{definition}

The RCP yields to the following results (see \cite{FT2003}, Theorems (10.7)
and (10.10)):

\begin{theorem}
\label{T1}
Let $G$ be a finite group acting on $\Lambda$ such that $\ker \tau$ has 
the Rotating Circle Property. 
Then local minima of $\action_\omega^G$ in $\Lambda^G$ do not
have any interior collision.
\end{theorem}
\begin{theorem}
\label{T2}
Let $G$ be a finite group acting on $\Lambda$ such that every maximal
$\T$-isotropy subgroup either has the RCP or acts trivially on
$\nn$, then any local minimizer of $\action_\omega^G$ in $\Lambda^G$
yields a collision-free periodic solution of the Newton 
equation~\ref{newton:eq}.
\end{theorem}

\begin{proof}[Proof of Theorem~\ref{MT2}]
We assume that $x(t)$ is a $G$-equivariant local minimizer.  If $\alpha \geq
2$ then $x(t)$ is collisionless, since the action level of colliding
trajectories is infinite.  When $\ker \tau$ is not trivial, 
Proposition~\ref{nr:ntc} 
ensures that  $x(t)$ is collisionless.  From now on we shall
assume that $\ker \tau$ is trivial and therefore $G$ is isomorphic to its
image  $\tau (G) \subset O(2)$.  On the other hand, Theorem~\ref{T1} says that
there are not interior collisions. This concludes the proof in the cyclic
case.

Let us now assume, by the sake of contradiction that there is a collision
instant at the boundary time $t=0$ and let $G_0$ be its $\mathbb{T}$-isotropy
subgroup, which is a group of order two (since $G \cong \tau (G)$) generated
by $g_0$. If $\sigma(g_0)$ acts trivially on the index set then 
Theorem~\ref{T2} leads to a contradiction.  Without loss of generality we can then
suppose that $\sigma(g_0)=(1,2)$.  As to $\rho(g_0)$, three cases are
possible: either $\rho(g_0)$ is the identity, the antipodal map or a
reflection with respect to a line $l$. In the first two cases the rotating
circle property holds and so we are left with the latter case.

Let $\kk\subset \nn = \{1,2,3\}$
be the colliding cluster;
it is easy to see that either $\kk=\{1,2\}$ 
or $\kk = \nn$. Let 
$q(t)=x(t)-x_0(t)$,
where $x_0(t)$ is the center of mass of the bodies in $\mathbf{k}$, and 
let $q^\lambda$ be defined by
\[
q^\lambda(t)=\lambda^{-2/(2+\alpha)}q(\lambda t).
\]
As shown in Section 7 of \cite{FT2003}, there are sequences $\lambda_n
\rightarrow 0$ such that $q^{\lambda_n}(t)$ converges to a $g_0$-equivariant
parabolic minimizing collision trajectory,
and this contradicts Theorem~\ref{main_thm}. 
The proof is now complete.
\end{proof}
 
\nocite{hsiang1994,hsiangstraume,hsiang2004,chenICM,or,%
montgomery,montgomery_contmath,%
sbano,marchal,chendesol,%
MR95k:58135,serter2,serter1,amco,bara91,cotizelati,%
amco94,bessi,%
argate,moeckel88,montgomery96,marchal_book,bruno}

\def\cfudot#1{\ifmmode\setbox7\hbox{$\accent"5E#1$}\else
  \setbox7\hbox{\accent"5E#1}\penalty 10000\relax\fi\raise 1\ht7
  \hbox{\raise.1ex\hbox to 1\wd7{\hss.\hss}}\penalty 10000 \hskip-1\wd7\penalty
  10000\box7} \def\cprime{$'$} \def\cprime{$'$} \def\cprime{$'$}
  \def\cprime{$'$} \def\cprime{$'$} \def\cprime{$'$} \def\cprime{$'$}
  \def\cprime{$'$} \def\cprime{$'$} \def\cprime{$'$} \def\cprime{$'$}
  \def\cprime{$'$} \def\cprime{$'$} \def\cprime{$'$} \def\cprime{$'$}
  \def\cprime{$'$} \def\polhk#1{\setbox0=\hbox{#1}{\ooalign{\hidewidth
  \lower1.5ex\hbox{`}\hidewidth\crcr\unhbox0}}}
  \def\polhk#1{\setbox0=\hbox{#1}{\ooalign{\hidewidth
  \lower1.5ex\hbox{`}\hidewidth\crcr\unhbox0}}} \def\cprime{$'$}
  \def\cprime{$'$} \def\cprime{$'$}


\appendix
\section*{Glossary}
Here  we provide a list
of terms together with the page at which the term
is introduced or described the first time.
\begin{description}
\item[Bound to collisions:] definition \ref{BTC} at page \pageref{BTC}.
\item[Coercive:] definition \ref{defi:Gcoervice} at page \pageref{defi:Gcoervice}. See also page \pageref{ap:coercive}.
\item[Cyclic, Brake or Dihedral type:] definition \ref{action-types} at page \pageref{action-types}.
\item[Fully uncoercive:] definition \ref{def:fullyunc} at page \pageref{def:fullyunc}.
\item[Homographic:] definition \ref{homG} at page \pageref{homG}.
\item[Interior or boundary collision:] definition \ref{defi:interiorboundarycollision} at page \pageref{defi:interiorboundarycollision}.
\item[Redundant:]  definition \ref{defi:redundant} at page \pageref{defi:redundant}.
\item[Rotating circle property:] definition \ref{rcp} at page \pageref{rcp}.
\item[$T$-isotropy:] definition \ref{T-isotropy} at page \pageref{T-isotropy}.
\item[Trivial core:] page \pageref{defi:core}.
\item[Type R:] definition \ref{def:typeR} at page \pageref{def:typeR}.
\end{description}

\end{document}